\documentclass[a4paper,12pt,reqno]{amsart}
\usepackage{amsmath}
\usepackage{amssymb}
\usepackage{amsthm}

\usepackage{amscd}
\usepackage{amsfonts}
\usepackage{setspace}
\usepackage{version}

\title[Partition function of zeta]{On the partition function of the Riemann zeta function, and the Fyodorov--Hiary--Keating conjecture}
\author{Adam J Harper}
\address{Mathematics Institute, Zeeman Building, University of Warwick, Coventry CV4 7AL, England}
\email{A.Harper@warwick.ac.uk}
\date{13th June 2019}

\numberwithin{equation}{section}
\theoremstyle{plain}

\onehalfspace
\newcommand{\N}{\mathbb{N}}
\newcommand{\R}{\mathbb{R}}
\newcommand{\E}{\mathbb{E}}
\newcommand{\p}{\mathbb{P}}
\newcommand{\Z}{\mathbb{Z}}

\newtheorem{thm1}{Theorem}
\newtheorem{thm2}[thm1]{Theorem}
\newtheorem{cor1}{Corollary}
\newtheorem{cor2}[cor1]{Corollary}

\newtheorem{lem1}{Lemma}
\newtheorem{lem2}[lem1]{Lemma}
\newtheorem{zetares1}{Zeta Function Result}
\newtheorem{lem3}[lem1]{Lemma}
\newtheorem{approxres1}{Approximation Result}

\newtheorem{keyprop1}{Key Proposition}
\newtheorem{keyprop2}[keyprop1]{Key Proposition}
\newtheorem{lem4}[lem1]{Lemma}

\newtheorem{lem5}[lem1]{Lemma}
\newtheorem{lem6}[lem1]{Lemma}
\newtheorem{probres1}{Probability Result}

\newtheorem{keyprop3}[keyprop1]{Key Proposition}
\newtheorem{keyprop4}[keyprop1]{Key Proposition}

\newtheorem{lem7}[lem1]{Lemma}

\begin{document}

\maketitle

\begin{abstract}
We investigate the ``partition function'' integrals $\int_{-1/2}^{1/2} |\zeta(1/2 + it + ih)|^2 dh$ for the critical exponent 2, and the local maxima $\max_{|h| \leq 1/2} |\zeta(1/2 + it + ih)|$, as $T \leq t \leq 2T$ varies. In particular, we prove that for $(1+o(1))T$ values of $T \leq t \leq 2T$ we have $\max_{|h| \leq 1/2} \log|\zeta(1/2+it+ih)| \leq \log\log T - (3/4 + o(1))\log\log\log T$, matching for the first time with both the leading and second order terms predicted by a conjecture of Fyodorov, Hiary and Keating.

The proofs work by approximating the zeta function in mean square by the product of a Dirichlet polynomial over smooth numbers and one over rough numbers. They then apply ideas and results from corresponding random model problems to compute averages of this product, under size restrictions on the smooth part that hold for most $T \leq t \leq 2T$ (but reduce the size of the averages). There are connections with the study of critical multiplicative chaos. Unlike in some previous work, our arguments never shift away from the critical line by more than a tiny amount $1/\log T$, and they don't require explicit calculations of Fourier transforms of Dirichlet polynomials.
\end{abstract}

\section{Introduction}
In this paper, our goal is to give upper bounds for different quantities involving the Riemann zeta function in short intervals on the critical line, that are sharp or close to sharp. As we shall explain, the original motivation for studying some of these quantities comes from probability or statistical mechanics, as do many of the proof ideas, but they all admit straightforward number theoretic definitions.

\vspace{12pt}
The first object we shall examine is the family of short integrals
$$ \int_{-1/2}^{1/2} |\zeta(1/2 + it + ih)|^2 dh = \int_{-1/2}^{1/2} e^{2\log|\zeta(1/2 + it + ih)|} dh , $$
where $T \leq t \leq 2T$, say. In the language of statistical mechanics, we can regard such an integral as something like a {\em partition function} corresponding to the {\em Hamiltonian} $\log|\zeta(1/2 + it + ih)|$ (where $t$ is given and $h$ varies). These integrals provide a certain measure of the average size of the zeta function in a short interval around the imaginary part $t$.

The classical second moment estimate for the zeta function (see e.g. Chapter 15 of Ivi\'c~\cite{ivic}) immediately implies that if we average the partition function by integrating over $t$, we have
$$ \frac{1}{T} \int_{T}^{2T} \Biggl( \int_{-1/2}^{1/2} |\zeta(1/2 + it + ih)|^2 dh \Biggr) dt \sim \log T \;\;\;\;\; \text{as} \; T \rightarrow \infty . $$
It is then natural to wonder whether, for most $T \leq t \leq 2T$, the partition function $\int_{-1/2}^{1/2} |\zeta(1/2 + it + ih)|^2 dh$ will have size around $\log T$. Our first theorem and corollary show this is not quite the case.

\begin{thm1}
Uniformly for all large $T$ and all $0 \leq q \leq 1$, we have
$$ \frac{1}{T} \int_{T}^{2T} \left( \int_{-1/2}^{1/2} |\zeta(1/2 + it+ih)|^2 dh \right)^q dt \ll \left( \frac{\log T}{1 + (1-q)\sqrt{\log\log T}} \right)^q . $$
\end{thm1}

\begin{cor1}
For all large $T$ and all $\lambda \geq 2$, we have
$$ \frac{1}{T} \text{meas}\{T \leq t \leq 2T : \int_{-1/2}^{1/2} |\zeta(1/2 + it+ih)|^2 dh \geq \lambda \frac{\log T}{\sqrt{\log\log T}} \} \ll \frac{\min\{\log \lambda, \sqrt{\log\log T}\}}{\lambda} . $$
\end{cor1}

\begin{proof}[Proof of Corollary 1]
If $\lambda \geq e^{\sqrt{\log\log T}}$, this follows immediately by applying Markov's inequality to the second moment estimate for zeta (or to Theorem 1 with $q=1$).

If $2 \leq \lambda \leq e^{\sqrt{\log\log T}}$, it follows by applying Markov's inequality to Theorem 1, with the choice $q = 1 - \frac{1}{2\log\lambda}$.
\end{proof}

Thus for most $T \leq t \leq 2T$, the partition function actually satisfies $\int_{-1/2}^{1/2} |\zeta(1/2 + it + ih)|^2 dh \ll \frac{\log T}{\sqrt{\log\log T}}$. The shape of the bound in Theorem 1 may look rather contrived. However, as we shall indicate when we discuss the proof, there is a good and interesting reason for the bound to take this form and it is likely to be best possible.

Fyodorov and Keating~\cite{fyodkeat} and Fyodorov, Hiary and Keating~\cite{fyodhiarykeat} considered the partition function integrals $\int_{-1/2}^{1/2} |\zeta(1/2 + it+ih)|^2 dh$, and the more general version $\int_{-1/2}^{1/2} |\zeta(1/2 + it+ih)|^{2\beta} dh$, as part of their investigations into the local maximum $\max_{|h| \leq 1/2} |\zeta(1/2 + it + ih)|$ that we shall discuss below. Their idea was that if one could compute all the ``moments of moments'' $\frac{1}{T} \int_{T}^{2T} \left( \int_{-1/2}^{1/2} |\zeta(1/2 + it+ih)|^{2\beta} dh \right)^q dt$ for varying $q$ then one would obtain strong distributional information about $\int_{-1/2}^{1/2} |\zeta(1/2 + it+ih)|^{2\beta} dh$, and if one could do this for varying $\beta$ one would obtain information about the local maximum. They calculated the analogous ``moments of moments'' for characteristic polynomials of random unitary matrices, for $q \in \N$ (so they could expand out the bracket) and $\beta < 1/\sqrt{q}$, using results on Toeplitz determinants. More recently, Bailey and Keating~\cite{baileykeat} extended the random matrix calculations for $q \in \N$ and $\beta \in \N$. However, this amount of information is not sufficient to rigorously draw conclusions about the local maximum (but very strong results about this in the random matrix setting are now known, see the papers of Arguin, Belius and Bourgade~\cite{abbrandmat}, of Paquette and Zeitouni~\cite{paqzeit}, and of Chhaibi, Madaule and Najnudel~\cite{chmadnaj} for increasingly precise theorems).

For the zeta function itself, Arguin, Ouimet and Radziwi\l\l~\cite{argouiradz} considered the integrals $\int_{-1/2}^{1/2} |\zeta(1/2 + it+ih)|^{2\beta} dh$, and in fact more general integrals over different interval lengths, for all fixed $\beta > 0$. They didn't quite consider the ``moments of moments'', but determined the size of these integrals up to factors of size $\log^{\epsilon}T$ for all $T \leq t \leq 2T$ apart from a set of size $o(T)$ (assuming the Riemann Hypothesis in some cases). These proofs exploit various interesting ideas, notably about the ``global'' moments of zeta and about the branching structure in its local maxima. But in the case of upper bounds when $\beta = 1$, they don't imply more about $\int_{-1/2}^{1/2} |\zeta(1/2 + it+ih)|^{2} dh$ than one gets by applying Markov's inequality to the classical second moment estimate.

This case where the exponent is 2 has a particular significance (when the integral is over an interval of fixed non-zero length) for the following reason: the dominant contribution to the $2\beta$-th moment of zeta comes from those $t$ with $|\zeta(1/2+it)| \approx \log^{\beta}T$ (see e.g. the introduction to Soundararajan's paper~\cite{soundmoments}). And for fixed $\beta < 1$ one expects to find many $t$ with zeta this large in each bounded interval, whilst for fixed $\beta > 1$ one expects to find no such $t$ in most intervals of bounded length. This explains the so-called {\em freezing transition} in the integrals $\int_{-1/2}^{1/2} |\zeta(1/2 + it+ih)|^{2\beta} dh$, conjectured by Fyodorov--Keating~\cite{fyodkeat} and Fyodorov--Hiary--Keating~\cite{fyodhiarykeat}, and proved up to factors $\log^{\epsilon}T$ by Arguin, Ouimet and Radziwi\l\l~\cite{argouiradz}. The transitional case where $\beta = 1$ is thus rather delicate, and the one which really captures information about the typical behaviour of the local maximum $\max_{|h| \leq 1/2} |\zeta(1/2 + it + ih)|$. Theorem 1 supplies us with some (probably) best possible information about this case.

\vspace{12pt}
Our second object of study will be the short interval maximum $\max_{|h| \leq 1/2} |\zeta(1/2 + it + ih)|$ that we already mentioned. This is the subject of a very precise conjecture by Fyodorov and Keating~\cite{fyodkeat} and Fyodorov, Hiary and Keating~\cite{fyodhiarykeat}, namely that for any real function $g(T)$ that tends to infinity with $T$, we should have
$$ \frac{1}{T} \text{meas}\{T \leq t \leq 2T : \bigl|\max_{|h| \leq 1/2} \log|\zeta(1/2+it+ih)| - (\log\log T - (3/4)\log\log\log T)\bigr| \leq g(T) \} \rightarrow 1 $$
as $T \rightarrow \infty$. In fact, their conjecture is an even more precise statement about the distribution of the difference $\max_{|h| \leq 1/2} \log|\zeta(1/2+it+ih)| - (\log\log T - (3/4)\log\log\log T)$. 

In the direction of this conjecture, using the second moment estimate for the Riemann zeta function and a Sobolev--Gallagher type argument (implying roughly that a large value of the zeta function should usually persist over an interval of length $\asymp 1/\log T$), it isn't too hard to show that for any $g(T)$ tending to infinity we have
$$ \frac{1}{T} \text{meas}\{T \leq t \leq 2T : \max_{|h| \leq 1/2} \log|\zeta(1/2+it+ih)| \leq \log\log T + g(T) \} \rightarrow 1 \;\;\; \text{as} \; T \rightarrow \infty . $$
See the papers of Arguin--Belius--Bourgade--Radziwi\l\l--Soundararajan~\cite{abbrs} and of Najnudel~\cite{najnudel}, as well as the end of the introduction to the preprint~\cite{harperlcz}. Here again we see the significance of the second moment of the zeta function when studying the local maximum. It is also known, but much harder (see later for some discussion of the proofs), that for any fixed $\epsilon > 0$ we have
$$ \frac{1}{T} \text{meas}\{T \leq t \leq 2T : \max_{|h| \leq 1/2} \log|\zeta(1/2+it+ih)| \geq (1-\epsilon)\log\log T \} \rightarrow 1 \;\;\; \text{as} \; T \rightarrow \infty . $$
This was proved by Najnudel~\cite{najnudel} assuming the truth of the Riemann Hypothesis, and by Arguin--Belius--Bourgade--Radziwi\l\l--Soundararajan~\cite{abbrs} unconditionally.

Here we look to improve the upper bound. In particular, as described by Fyodorov--Keating~\cite{fyodkeat} and Fyodorov--Hiary--Keating~\cite{fyodhiarykeat} (see also the author's Bourbaki survey~\cite{harperbourbaki} for a gentle discussion), the second order term $- (3/4)\log\log\log T$ in the conjecture would be a manifestation of some dependence (or non-trivial correlations) between the values of $\zeta(1/2 + it + ih)$ at nearby $h$. If these values behaved independently at spacings of about $1/\log T$, one would instead expect the maximum to usually be around $\log\log T - (1/4)\log\log\log T$.

\begin{thm2}
Uniformly for all large $T$ and $0 \leq U \leq \log\log T$, we have
\begin{eqnarray}
&& \frac{1}{T} \text{meas}\{T \leq t \leq 2T : \max_{|h| \leq 1/2} |\zeta(1/2+it+ih)| \geq \frac{e^{U} \log T}{(\log\log T)^{3/4}} \} \nonumber \\
& \ll & e^{-2U} (\log\log\log T + U) (\log\log\log T)^2 . \nonumber
\end{eqnarray}
\end{thm2}

\begin{cor2}
For any real function $g(T)$ tending to infinity with $T$, we have
$$ \max_{|h| \leq 1/2} \log|\zeta(1/2+it+ih)| \leq \log\log T - (3/4)\log\log\log T + (3/2)\log\log\log\log T + g(T) $$
for a set of $T \leq t \leq 2T$ with measure $(1+o(1))T$.
\end{cor2}

Corollary 2 matches as an upper bound the first two terms in the Fyodorov--Hiary--Keating conjecture, showing that the alternative ``independence conjecture'' $\log\log T - (1/4)\log\log\log T$ for the typical size of the short interval maximum cannot be correct. Unfortunately it still (presumably) isn't quite sharp thanks to the third order term $(3/2)\log\log\log\log T$, needed to overcome the powers of $\log\log\log T$ in Theorem 2. The sharpest distributional form of the Fyodorov--Hiary--Keating conjecture~\cite{fyodkeat, fyodhiarykeat} suggests that the correct Theorem 2 upper bound is of the shape $Ue^{-2U}$, at least for large fixed $U$. The primary source, although not the only source, of the extraneous $\log\log\log T$ powers in Theorem 2 is a crude treatment of large primes contributing to the zeta function, which effectively decreases the threshold $\frac{e^{U} \log T}{(\log\log T)^{3/4}}$ in the course of the proof (by forcing a parameter $V$ to be overly large, see below).

\vspace{12pt}
{\em Remark.} In work in preparation, Arguin, Bourgade, Radziwi\l\l \; and Soundararajan have independently obtained a proof of Corollary 2. The methods of their proof are rather different than those developed here.

\subsection{Ideas from the proofs}
Next we describe our strategy for proving the theorems. Later we shall make some comparisons with the previous work on the Fyodorov--Hiary--Keating conjecture. Let $(f(p))_{p \; \text{prime}}$ be a sequence of independent Steinhaus random variables, i.e. independent random variables distributed uniformly on the complex unit circle. Then define the random Euler product $F(s) := \prod_{p \leq x} (1 - \frac{f(p)}{p^s})^{-1}$, where $x$ is a large parameter. In the author's work~\cite{harperrmflow} on low moments of random multiplicative functions, a key issue in the proofs was to show that, uniformly for $0 \leq q \leq 1$, we have
$$ \E\Biggl(\int_{-1/2}^{1/2} |F(1/2+ih)|^2 dh \Biggr)^{q} \asymp \left(\frac{\log x}{1 + (1-q)\sqrt{\log\log x}}\right)^q . $$
It is a well known principle in analytic number theory, and in particular in work on the value distribution of $L$-functions and on the Fyodorov--Hiary--Keating conjecture, that as $t$ varies over a long interval the numbers $(p^{-it})_{p \; \text{prime}}$ should ``behave like'' a sequence of independent Steinhaus random variables. It is also a well known principle that, in many respects, the statistical behaviour of $\zeta(1/2+it+ih)$ should resemble that of an Euler product $\prod_{p \leq t} (1 - \frac{1}{p^{1/2+it+ih}})^{-1}$. From this perspective, Theorem 1 is simply the derandomised version of the (upper bound) result we already have in the random case.

\vspace{12pt}
In order to accomplish this derandomisation rigorously, we proceed in several steps. Firstly, using the approximate functional equation we can upper bound $|\zeta(1/2+it)|$ by something like $|\sum_{n \leq \sqrt{t/2\pi}} \frac{1}{n^{1/2+it}}|$, up to very small error terms. In fact, for the proof of Theorem 1 it would suffice to use the simpler Hardy--Littlewood approximation to replace $\zeta(1/2+it)$ by $\sum_{n \leq T} \frac{1}{n^{1/2+it}}$, but this would not work for the proof of Theorem 2 and we prefer to give a unified treatment so far as possible.

Next, we want to replace $\sum_{n \leq \sqrt{t/2\pi}} \frac{1}{n^{1/2+it}}$ by something more like an Euler product, to bring us closer to the random setting. We are helped by the fact that the saving $1 + (1-q)\sqrt{\log\log T}$ we are seeking, as compared with a ``trivial'' bound using the second moment of zeta, is a very slowly growing function of $T$, so rather than showing that the whole sum $\sum_{n \leq \sqrt{t/2\pi}} \frac{1}{n^{1/2+it}}$ behaves like a product it will suffice to show that a not too small piece of it does. By a fairly simple argument using the mean value theorem for Dirichlet polynomials (see Lemma \ref{polyeditlem}, below), we will show that up to acceptable error terms we usually have
$$ \sum_{n \leq \sqrt{t/2\pi}} \frac{1}{n^{1/2+it}} \approx \sum_{\substack{m \leq T^{\epsilon}, \\ m \; \text{is} \; P \; \text{smooth}}} \frac{1}{m^{1/2+it}} \sum_{\substack{n \leq T^{1/2 - 2\epsilon}, \\ n \; \text{is} \; P \; \text{rough}}} \frac{1}{n^{1/2+it}} . $$
Recall that a number is said to be {\em $P$-smooth} if all of its prime factors are $\leq P$, and to be {\em $P$-rough} if all of its prime factors are $> P$. For the proof of Theorem 1 there is considerable flexibility in the choice of $P, \epsilon$ (there is less for Theorem 2), but we take $P = T^{1/(\log\log T)^8}$ and $\epsilon = \frac{1}{(\log\log T)^2}$. Because the density of the $P$-smooth numbers that are $> T^{\epsilon}$ is small, the sum over $m$ will turn out to behave much like the full product $\prod_{p \leq P} (1 - \frac{1}{p^{1/2+it}})^{-1} = \sum_{\substack{m = 1, \\ m \; \text{is} \; P \; \text{smooth}}}^{\infty} \frac{1}{m^{1/2+it}}$. Some other nice features here (which are probably not essential, but very convenient) are the fact that the ranges of summation over $m$ and $n$ do not interact with one another, and the fact that the maximum size of product $mn$ appearing is $T^{\epsilon} T^{1/2-2\epsilon} = T^{1/2-\epsilon}$, in other words we have reduced the total length of summation from about $T^{1/2}$ to $T^{1/2 - \epsilon}$. We can make all these edits to the sum because the upper bound we are trying to prove is fairly large (roughly $\log T$), so we can tolerate moderately large error terms in our approximations.

Now we need to look inside the proofs from the random case. To upper bound $\E(\int_{-1/2}^{1/2} |F(1/2+ih)|^2 dh )^{q}$, one considers an event $\mathcal{G}$ that neither $|F(1/2+ih)|$ nor its partial products are too large (depending on $q$ and on the length of the product) for any $|h| \leq 1/2$. Roughly speaking, one shows that the event $\mathcal{G}$ will occur with very high probability, and also that $\E \textbf{1}_{\mathcal{G}} \int_{-1/2}^{1/2} |F(1/2+ih)|^2 dh \ll \frac{\log x}{1 + (1-q)\sqrt{\log\log x}}$. These two estimates may be combined in a recursive procedure, varying $q$, to prove the desired upper bound for the $q$-th moment. The introduction to the author's paper~\cite{harperrmflow} contains much more discussion of the origins and motivation of this argument. Let us note, however, that the key force at work is the close relationship between $\int_{-1/2}^{1/2} |F(1/2+ih)|^2 dh$ and a probabilistic object called {\em critical multiplicative chaos}. It is {\em critical} because the exponent 2 of the product is the one at which the integral starts to be dominated by very rare events, as discussed earlier. The denominator $\sqrt{\log\log x}$ appearing in our results reflects a non-trivial renormalisation that one can perform when studying critical chaos (but not in the subcritical case of exponent $< 2$). See also the papers of Saksman and Webb~\cite{saksmanwebb, saksmanwebb2} for further work on connections between the Riemann zeta function and multiplicative chaos.

In our deterministic setting, we consider an analogous event $\mathcal{G}_t$ that neither $\prod_{p \leq P} (1 - \frac{1}{p^{1/2+it+ih}})^{-1}$ nor its partial products are too large, for any $|h| \leq 1/2$. The proof that the random event $\mathcal{G}$ has high probability was a relatively simple argument using the union bound and exponential moment estimates, so we can duplicate this using the union bound and high moment estimates for Dirichlet polynomials over primes, provided the parameter $P$ isn't too large compared with $T$. The estimation of $\E \textbf{1}_{\mathcal{G}} \int_{-1/2}^{1/2} |F(1/2+ih)|^2 dh$ is probabilistically much more subtle, but fortunately in the analogous calculation of $\frac{1}{T} \int_{T}^{2T} \textbf{1}_{\mathcal{G}_t} \int_{-1/2}^{1/2}  |\sum_{\substack{m \leq T^{\epsilon}, \\ m \; \text{is} \; P \; \text{smooth}}} \frac{1}{m^{1/2+it+ih}} \sum_{\substack{n \leq T^{1/2 - 2\epsilon}, \\ n \; \text{is} \; P \; \text{rough}}} \frac{1}{n^{1/2+it+ih}}|^2 dh dt$ we can piggyback on that existing probabilistic result. Roughly speaking, we approximate the indicator function $\textbf{1}_{\mathcal{G}_t}$ by a product of smooth functions evaluated at various Dirichlet polynomials over primes of size at most $P$. We then Taylor expand each smooth function sufficiently far that the contribution from the remainder is small (expanding to order $(\log\log T)^{O(1)}$ is sufficient). The remaining integrand is simply a sum of products of Dirichlet polynomials, and provided $P$ isn't too large compared with $T$, mean value estimates show this is the same (up to very small error terms) as the corresponding expectation with $p^{-it}, m^{-it}, n^{-it}$ replaced everywhere by the values $f(p), f(m), f(n)$ of a Steinhaus random multiplicative function\footnote{A {\em Steinhaus random multiplicative function} $f(n)$ is the random totally multiplicative function obtained by extending the sequence $(f(p))_{p \; \text{prime}}$ of independent Steinhaus random variables. In other words, we define $f(n) := \prod_{p^a ||n} f(p)^a$ for all natural numbers $n$, where $p^a ||n$ means that $p^a$ is the highest power of the prime $p$ that divides $n$.}. Having passed to this random case, we can undo the previous steps and also replace the random sum $\sum_{\substack{m \leq T^{\epsilon}, \\ m \; \text{is} \; P \; \text{smooth}}} \frac{f(m)}{m^{1/2+ih}}$ by the product $\prod_{p \leq P} (1 - \frac{f(p)}{p^{1/2+ih}})^{-1}$, which it very closely approximates in mean square by choice of $\epsilon$. Finally, because the random version of the event $\mathcal{G}_t$ only depends on the $(f(p))_{p \leq P}$, by independence we can simply replace $|\sum_{\substack{n \leq T^{1/2 - 2\epsilon}, \\ n \; \text{is} \; P \; \text{rough}}} \frac{f(n)}{n^{1/2+ih}}|^2$ by its expectation, which is $\asymp (\log T)/\log P$. The remaining integral is (more or less) exactly of the form $\E \textbf{1}_{\mathcal{G}} \int_{-1/2}^{1/2} |F(1/2+ih)|^2 dh$, with $x$ replaced by $P$, and deploying the existing estimate for that completes the proof.

\vspace{12pt}
For Theorem 2, much is the same but a few additional ingredients are required. We must pass from looking at $\max_{|h| \leq 1/2} |\zeta(1/2 + it + ih)|$, which is the value of the zeta function at some unknown point with imaginary part in the interval $[t-1/2,t+1/2]$, to a version that is easier to average and connect with the proof of Theorem 1. Many approaches could probably be effective for this, but we take a simple one (apparently not used for these problems before) and use Cauchy's Integral Formula to replace the maximum by an integral around a small rectangle (of side length $\asymp 1/\log T$) about the point where the maximum is attained. At this step the position of this point remains unknown, but later in the argument we can deal with this by extending the lengths of the vertical integrals and summing over all possible positions of the horizontal integrals. Although this appears wasteful (it corresponds to applying the union bound on the probabilistic side), in fact it is efficient because large values of the zeta function are rare and these make the dominant contribution to the integrals being summed.

Having passed to integrals of zeta, we can replace the zeta function by an integrand roughly of the form $\sum_{\substack{m \leq T^{\epsilon}, \\ m \; \text{is} \; P \; \text{smooth}}} \frac{1}{m^{1/2+it+ih}} \sum_{\substack{n \leq T^{1/2 - 2\epsilon}, \\ n \; \text{is} \; P \; \text{rough}}} \frac{1}{n^{1/2+it+ih}}$, as in the proof of Theorem 1. Since the event $\mathcal{G}_t$ occurs for most $T \leq t \leq 2T$, we can also restrict to considering those $t$ for which it holds. In fact we restrict to those $t$ for which a slightly modified event $\tilde{\mathcal{G}}_t$ holds, whose definition in particular depends on $U$, so that we can obtain a measure bound depending on $U$.

The next issue is that a priori we are looking for values of the zeta function, or equivalently of our small rectangular integral, that are $\geq \frac{e^{U} \log T}{(\log\log T)^{3/4}}$, but we want to do most of our work with the smooth part $\sum_{\substack{m \leq T^{\epsilon}, \\ m \; \text{is} \; P \; \text{smooth}}} \frac{1}{m^{1/2+it+ih}}$ of the integrand (as that is the part which closely resembles an Euler product). So we need to know that if the integral is large, then generally the smooth part must be remarkably large as well. It certainly seems reasonable to expect this, as the mean square average of the smooth part is $\sum_{\substack{m \leq T^{\epsilon}, \\ m \; \text{is} \; P \; \text{smooth}}} \frac{1}{m} \asymp \log P$, whereas the mean square average of the rough sum is the much smaller quantity $\sum_{\substack{n \leq T^{1/2 - 2\epsilon}, \\ n \; \text{is} \; P \; \text{rough}}} \frac{1}{n} \asymp \frac{\log T}{\log P}$. So if the sum over $m$ were not large it would be very difficult for the sum over $n$ to compensate and produce a very large product. To make this rigorous, we look at the {\em fourth} moment of $\sum_{\substack{m \leq T^{\epsilon}, \\ m \; \text{is} \; P \; \text{smooth}}} \frac{1}{m^{1/2+it}} \sum_{\substack{n \leq T^{1/2 - 2\epsilon}, \\ n \; \text{is} \; P \; \text{rough}}} \frac{1}{n^{1/2+it}}$ restricted to $t$ for which $|\sum_{\substack{m \leq T^{\epsilon}, \\ m \; \text{is} \; P \; \text{smooth}}} \frac{1}{m^{1/2+it}}| \leq \frac{\log P}{V}$, see Lemma \ref{fourthmomentlem} below. (This is fairly easy, but it is the application of the mean value theorem for Dirichlet polynomials at this point, and only here, where it is important that the maximum size of $mn$ is no more than $\sqrt{T}$). The bound obtained has a quadratic ``saving'' $V^{-2}$, along with a ``loss'' of the form $\frac{\log T}{\log P} = (\log\log T)^8$ coming from the rough Dirichlet polynomial over $n$. So if $V$ is a suitable power of $\log\log T$ (or, when $U$ is large, a power of $\log\log T$ times $e^{-U}$), the overall contribution is sufficiently small, and we can restrict thereafter to points where $|\sum_{\substack{m \leq T^{\epsilon}, \\ m \; \text{is} \; P \; \text{smooth}}} \frac{1}{m^{1/2+it}}| > \frac{\log P}{V}$. Note that this choice of $V$ is a key reason that our bound in Theorem 2 isn't ultimately sharp.

At this point our situation is close to Theorem 1, except that in our integrals
$$ \frac{1}{T} \int_{T}^{2T} \textbf{1}_{\tilde{\mathcal{G}}_t} \int_{-1/2}^{1/2} \textbf{1}_{|\sum_{\substack{m \leq T^{\epsilon}, \\ P \; \text{smooth}}} \frac{1}{m^{1/2+it+ih}}| > (\log P)/V} |\sum_{\substack{m \leq T^{\epsilon}, \\ P \; \text{smooth}}} \frac{1}{m^{1/2+it+ih}} \sum_{\substack{n \leq T^{1/2 - 2\epsilon}, \\ P \; \text{rough}}} \frac{1}{n^{1/2+it+ih}}|^2 dh dt $$
we don't just have the $\tilde{\mathcal{G}}_t$ restriction that $\prod_{p \leq P} (1 - \frac{1}{p^{1/2+it+ih}})^{-1}$ and its partial products aren't too large, but also the restriction that $|\sum_{\substack{m \leq T^{\epsilon}, \\ P \; \text{smooth}}} \frac{1}{m^{1/2+it+ih}}| > (\log P)/V$. With a little more work, see Lemma \ref{eulerprodapproxlem} below, this lower bound condition can be replaced by a condition roughly like $\prod_{p \leq P} |1 - \frac{1}{p^{1/2+it+ih}}|^{-1} \gg (\log P)/V$. Whereas the upper bound restriction from $\tilde{\mathcal{G}}_t$ reduces the size of the integral by around a factor $\frac{1}{\sqrt{\log\log P}}$, as in Theorem 1, the additional lower bound reduces it overall by a factor around $\frac{1}{(\log\log P)^{3/2}}$, see Probability Result 1 below. Since the integrand roughly corresponds to the square of the zeta function, and we chose parameters such that $\log\log P \asymp \log\log T$, this directly suggests the denominator $(\log\log T)^{3/4}$ in Theorem 2. To conclude rigorously we argue as in the proof of Theorem 1, using Taylor expansion and the mean value theorem for Dirichlet polynomials to replace the integral by the corresponding average of random multiplicative $f(n)$, and then applying Probability Result 1 and relevant machinery from the author's paper~\cite{harperrmflow}. Conceptually this part is fairly straightforward, although some technical work is required to handle a discretisation, and various size restrictions, in the definition of $\tilde{\mathcal{G}}_t$ and in the machinery from \cite{harperrmflow}.

\subsection{Comparisons and possible improvements}
A main difference between the proofs of Theorems 1 and 2, and previous work on the Fyodorov--Hiary--Keating conjecture, lies in the way we approximate the zeta function by Dirichlet polynomials. In our approximations we keep the term $\sum_{\substack{n \leq T^{1/2 - 2\epsilon}, \\ n \; \text{is} \; P \; \text{rough}}} \frac{1}{n^{1/2+it}}$ around. This is never shown to behave like an Euler product (indeed it won't really do so), it simply remains until the end where it is estimated in mean square. In the work of Najnudel~\cite{najnudel} and of Arguin--Belius--Bourgade--Radziwi\l\l--Soundararajan~\cite{abbrs}, they instead lower bound $\max_{|h| \leq 1/2} |\zeta(1/2+it+ih)|$ solely by the maximum of the exponential of a prime number sum (i.e. an Euler product, essentially), but to do this they shorten the sum/product and throw away the contribution from larger primes. (In Najnudel's work this shortening is explicit, and the Riemann Hypothesis is invoked to make it work for {\em every} $t$. For Arguin--Belius--Bourgade--Radziwi\l\l--Soundararajan, the shortening arises from shifting to look at $\zeta(\sigma+it+ih)$ with $\sigma$ slightly away from 1/2, and the Riemann Hypothesis is avoided because one only looks to do so for {\em most} $t$.) These procedures reduce the problem to distributional calculations for Dirichlet polynomials over primes, but the shortening lowers the size of maximum one can detect. In contrast, our arguments entail extra work to relate the desired lower bound condition on zeta to one for the smooth Dirichlet polynomial $\sum_{\substack{m \leq T^{\epsilon}, \\ P \; \text{smooth}}} \frac{1}{m^{1/2+it+ih}}$, but lose less because the large primes haven't been discarded. We also gain some simplifications from piggybacking our calculations on those already performed in the random case~\cite{harperrmflow}, rather than computing things again for a particular set of prime number sums. An important caveat is that if one were looking to prove lower bound counterparts of Theorem 2 and Corollary 2, it is not clear whether one could successfully pass between bounds for zeta and bounds for $\sum_{\substack{m \leq T^{\epsilon}, \\ P \; \text{smooth}}} \frac{1}{m^{1/2+it+ih}}$ in this fashion. See below for some further explanation of this.

It would be nice to improve the upper bounds in Theorem 2 and Corollary 2 so that they were really sharp. As already remarked, our biggest inefficiency at the moment comes from Lemma \ref{fourthmomentlem}, where we use a fourth moment argument to export the lower bound condition on zeta to a lower bound condition on $\sum_{\substack{m \leq T^{\epsilon}, \\ m \; \text{is} \; P \; \text{smooth}}} \frac{1}{m^{1/2+it}}$, incurring some loss $\asymp (\log T)/\log P$ from the fourth moment of the rough part $\sum_{\substack{n \leq T^{1/2 - 2\epsilon}, \\ n \; \text{is} \; P \; \text{rough}}} \frac{1}{n^{1/2+it}}$. An obvious modification is to approximate zeta by a product of more Dirichlet polynomials, say one over $P$-smooth numbers, one over numbers with all their prime factors between $P$ and $P'$, and one over $P'$-rough numbers, for suitable $P'$ much larger than $P$. Because $P'$ is larger than $P$, a fourth moment argument could let one impose a condition that the {\em product} of the first two Dirichlet polynomials must be large, with less loss (i.e. smaller $V$) than currently. And by Lemma \ref{fourthmomentlem} we already know that the $P$-smooth Dirichlet polynomial must be fairly large, so we would have quite good control on the size of the middle Dirichlet polynomial for the rest of the argument. The author finds it quite plausible that such an argument could improve our results (but has checked no details), but getting a bound that is sharp down to order one terms appears very challenging. In the first place one would probably need to approximate zeta by the product of a growing number of Dirichlet polynomials. We also note that if one wants such a sharp result, new difficulties appear on the more probabilistic sides of the argument as well. Note the substantial work~\cite{abbrandmat, paqzeit, chmadnaj} that was required to achieve this in the random matrix case, whilst for the Steinhaus random multiplicative model of the zeta function the analogue of Corollary 2 is currently only known up to the second order term~\cite{abh}.

Another goal is to obtain matching lower bounds in the theorems. As we mentioned, the lower bound $\E(\int_{-1/2}^{1/2} |F(1/2+ih)|^2 dh )^{q} \gg (\frac{\log x}{1 + (1-q)\sqrt{\log\log x}} )^q$ is already known for the random analogue of Theorem 1, and it is plausible that by combining the methods of this paper with that probabilistic proof one could obtain lower bounds for $\frac{1}{T} \int_{T}^{2T} ( \int_{-1/2}^{1/2} |\zeta(1/2 + it+ih)|^2 dh )^q dt$ as well. Although the lower bound argument is more involved than the upper bound, it would ultimately reduce to {\em upper} bounding averages of $|\zeta(1/2 + it+ih_1)|^2 |\zeta(1/2 + it+ih_2)|^2$ subject to growth constraints like $\mathcal{G}_t$, which seems achievable. But for a lower bound analogue of Theorem 2, one would not only need all this with an additional constraint that $\sum_{\substack{m \leq T^{\epsilon}, \\ P \; \text{smooth}}} \frac{1}{m^{1/2+it+ih}}$ should be large (which again might be achievable), but also to deduce from a lower bound for these restricted integrals that there must be a very large value of zeta on most intervals $[t-1/2,t+1/2]$. It is not so clear how to distinguish this from the case of several moderately large zeta values contributing to the integral, without more information about the rough Dirichlet polynomial $\sum_{\substack{n \leq T^{1/2 - 2\epsilon}, \\ n \; \text{is} \; P \; \text{rough}}} \frac{1}{n^{1/2+it+ih}}$ which is difficult to control.

\subsection{Organisational remarks}
In section 2, we collect a few results about Dirichlet polynomials and approximating functions that are needed to prove both Theorems 1 and 2. Section 3 contains the proof of Theorem 1, which is essentially a derandomised version of the upper bound part of Theorem 1 of Harper~\cite{harperrmflow}. In section 4 we present two further tools needed specifically for the proof of Theorem 2, namely a fourth moment estimate and a certain probabilistic estimate concerning Gaussian random walks (a version of the Ballot Theorem). Finally, in section 5 we deploy our tools to prove Theorem 2. Some of the derandomisation steps there are very close to those from the proof of Theorem 1, so not presented in detail.

\section{Tools and preliminary results}
A basic tool that we shall require is a standard kind of mean value estimate for Dirichlet polynomials.

\begin{lem1}\label{meanvaluelem}
For any $T, H > 0$, any $x \geq 1$, and any complex numbers $(a(n))_{n \leq x}, (b(n))_{n \leq x}$, we have
$$ \frac{1}{H} \int_{T}^{T + H} \left(\sum_{n \leq x} a(n) n^{-it} \right) \overline{\left( \sum_{m \leq x} b(m) m^{-it} \right)} dt = \sum_{n \leq x} a(n) \overline{b(n)} + O\Biggl(\frac{x}{H} \sqrt{\sum_{n \leq x} |a(n)|^2 \sum_{m \leq x} |b(m)|^2} \Biggr) . $$

In particular, for any $1 \leq x \leq H$ and any complex numbers $(a(n))_{n \leq x}$, we have
$$ \frac{1}{H} \int_{T}^{T + H} \left|\sum_{n \leq x} a(n) n^{-it} \right|^2 dt \ll \sum_{n \leq x} |a(n)|^2 . $$
\end{lem1}
The first statement follows by expanding out the sums and integral and using Hilbert's Inequality to bound the total contribution from all the terms with $m \neq n$. See chapter 7.5 of Montgomery~\cite{mont}, for example. Note that the first term on the right equals $\E (\sum_{n \leq x} a(n) f(n)) \overline{(\sum_{m \leq x} b(m) f(m))}$, where $f(n)$ is a Steinhaus random multiplicative function. The second statement is an immediate corollary of the first, on choosing the sequence $(b(n))_{n \leq x}$ to be the same as the sequence $(a(n))_{n \leq x}$.

Applying this mean value estimate in a fairly routine way (apart from a little fiddling about with divisor functions and squares of primes), we obtain the following upper bound for the even integer moments of various Dirichlet polynomials. We will use this to obtain large deviation estimates for Dirichlet polynomials, and also to control the contribution from remainder terms in Taylor expansions.

\begin{lem2}\label{highmomentslem}
Let $x \geq 1$, and let $(c(n))_{n \leq x}$ be any complex numbers. Let $\mathcal{P}$ be any finite set of primes, let $\mathcal{Q}$ be any (non-empty) set consisting of some elements of $\mathcal{P}$ and squares of elements of $\mathcal{P}$, and write $U := \max\{q \in \mathcal{Q}\}$ . Finally, let $Q(t) := \sum_{q \in \mathcal{Q}} \frac{a(q)}{q^{1/2+it}}$, where the $a(q)$ are any complex numbers.

Then for any natural number $k$ such that $x U^k < T$, we have
$$ \frac{1}{T} \int_{T}^{2T} \Biggl|\sum_{n \leq x} c(n) n^{-it} \Biggr|^{2} |Q(t)|^{2k} dt \ll \Biggl(\sum_{n \leq x} \tilde{d}(n) |c(n)|^2 \Biggr) \cdot (k !) \Biggl( 2 \sum_{q \in \mathcal{Q}} \frac{v_q |a(q)|^2}{q} \Biggr)^{k} , $$
where $\tilde{d}(n) := \sum_{d|n} \textbf{1}_{p|d \Rightarrow p \in \mathcal{P}}$, and $v_q$ is 1 if $q$ is a prime and 6 if $q$ is the square of a prime.

Also, for any natural number $k$ such that $U^{k} < T$, we have the sharper bound
$$ \frac{1}{T} \int_{T}^{2T} |Q(t)|^{2k} dt \ll (k !) \Biggl( \sum_{q \in \mathcal{Q}} \frac{v_q |a(q)|^2}{q} \Biggr)^{k} . $$
\end{lem2}

\begin{proof}[Proof of Lemma \ref{highmomentslem}]
For the first part, we rewrite the left hand side as
$$ \frac{1}{T} \int_{T}^{2T} \Biggl|\sum_{u \leq xU^k} u^{-it} \sum_{\substack{nm=u, \\ n \leq x, \\ m \leq U^k}} \frac{c(n) a(m)}{\sqrt{m}} \Biggr|^{2} dt , $$
where $a(m) := \sum_{q_{1}, ..., q_{k} \in \mathcal{Q}} \textbf{1}_{\prod q_{i} = m} \prod_{i} a(q_{i})$. Note in particular that $a(m)$ is supported on numbers $m$ having all their prime factors from $\mathcal{P}$. Applying the second part of Lemma \ref{meanvaluelem}, and using our condition that $x U^k < T$, we find this is all
$$ \ll \sum_{u \leq x U^{k}} \Biggl| \sum_{\substack{nm=u, \\ n \leq x, \\ m \leq U^k}} \frac{c(n) a(m)}{\sqrt{m}} \Biggr|^{2} . $$
Next, the Cauchy--Schwarz inequality applied to the sum over $n$ and $m$, and the sub-multiplicativity of the function $\tilde{d}(n)$, imply the above is
\begin{eqnarray}
\leq \sum_{u \leq x U^k} \tilde{d}(u) \sum_{\substack{nm=u, \\ n \leq x, \\ m \leq U^k}} \frac{|c(n)|^2 |a(m)|^2}{m} & \leq & \sum_{u \leq x U^k} \sum_{\substack{nm=u, \\ n \leq x, \\ m \leq U^k}} \frac{\tilde{d}(n) |c(n)|^2 \tilde{d}(m) |a(m)|^2}{m} \nonumber \\
& = & \left( \sum_{n \leq x} \tilde{d}(n) |c(n)|^2 \right) \left( \sum_{m \leq U^k} \frac{\tilde{d}(m) |a(m)|^2}{m} \right) . \nonumber
\end{eqnarray}

The Cauchy--Schwarz inequality and the sub-multiplicativity of $\tilde{d}(m)$ further imply that $\frac{\tilde{d}(m) |a(m)|^2}{m}$ is
$$ \leq \frac{\tilde{d}(m)}{m}(\sum_{q_1, ..., q_k \in \mathcal{Q}} \textbf{1}_{\prod q_{i} = m} ) (\sum_{\substack{q_1, ..., q_k \in \mathcal{Q}, \\ \prod q_i = m}} \prod_{i} |a(q_{i})|^2) \leq (\sum_{q_1, ..., q_k \in \mathcal{Q}} \textbf{1}_{\prod q_{i} = m} ) (\sum_{\substack{q_1, ..., q_k \in \mathcal{Q}, \\ \prod q_i = m}} \prod_{i} \frac{\tilde{d}(q_i) |a(q_{i})|^2}{q_{i}}) . $$
Now if $m$ has prime factorisation $m = \prod_{j=1}^{w} p_{j}^{l_j}$, where the $p_j$ are distinct primes, then the collection of all possible multisets of $k$ primes and prime squares with product $m$ is parametrised by tuples $(r_1, ..., r_w)$, where $0 \leq r_j \leq l_j$ denotes the number of {\em single} copies of $p_j$ in the multiset. (Thus $r_j$ must have the same parity as $l_j$, and there will be $(l_j - r_j)/2$ copies of $p_j^2$ in the multiset. In particular, the number of multisets is crudely at most $\prod_{j=1}^{w} l_j$.) For any such multiset, the number of tuples $q_1, ..., q_k \in \mathcal{Q}$ whose elements are those of the multiset (in some order) is very crudely at most $(k!)/\prod_{j=1}^{w} (r_j !) \leq (k!)/\prod_{j=1}^{w} 2^{r_j - 1} = (k!) \prod_{j=1}^{w} 2^{l_j - r_j} /\prod_{j=1}^{w} 2^{l_j - 1} \leq (k!) \prod_{j=1}^{w} 2^{l_j - r_j} /\prod_{j=1}^{w} l_j$. And both of the products here depend only on $m$ (not on the particular multiset), in particular $\prod_{j=1}^{w} 2^{l_j - r_j} = 4^{\#\{\text{squares in the multiset}\}} = 4^{\Omega(m) - k}$, where $\Omega(m) = \sum_{j=1}^{w} l_j$ is the total number of prime factors of $m$.

Putting everything together, recalling the definition of $v_q$ and that $\tilde{d}(q_i)$ equals 2 if $q_i \in \mathcal{Q}$ is prime and equals 3 if $q_i \in \mathcal{Q}$ is the square of a prime, we conclude that
$$ \frac{\tilde{d}(m) |a(m)|^2}{m} \leq (k!) 4^{\Omega(m) - k} \sum_{\substack{q_1, ..., q_k \in \mathcal{Q}, \\ \prod q_i = m}} \prod_{i} \frac{\tilde{d}(q_i) |a(q_{i})|^2}{q_{i}} = (k!) \sum_{\substack{q_1, ..., q_k \in \mathcal{Q}, \\ \prod q_i = m}} \prod_{i} \frac{2 v_{q_i} |a(q_{i})|^2}{q_{i}} . $$
The first part of Lemma \ref{highmomentslem} follows on summing over $m$.

In the special case where $x=1$ and $c(1)=1$, there is no need to apply the Cauchy--Schwarz inequality to the sum over $n$ and $m$ (it only has 1 term). Thus we don't pick up any factor $\tilde{d}(m)$ in the argument, and don't end up with a multiplier $2$ at the end. This second part is also Lemma 3 of Soundararajan~\cite{soundmoments} (under slightly more general conditions).
\end{proof}

\vspace{12pt}
In order to work with the Riemann zeta function, we shall use the approximate functional equation (more specifically, the symmetric form on the critical line) to replace it by Dirichlet polynomials of suitable length.

\begin{zetares1}[Approximate Functional Equation]
For all large real $t$, we have
$$ \zeta(1/2 + it) = \sum_{n \leq \sqrt{t/2\pi}} \frac{1}{n^{1/2+it}} + \chi(1/2 + it) \sum_{n \leq \sqrt{t/2\pi}} \frac{1}{n^{1/2-it}} + O(t^{-1/4}) , $$
where the function $\chi(\cdot)$ satisfies $|\chi(1/2+it)| = 1$ for all $t \in \R$.
\end{zetares1}

See e.g. Chapter 4.1 of Ivi\'c~\cite{ivic} for a proof of this.

Using our mean value estimate, we can (in an average sense) replace the Dirichlet polynomials supplied by the approximate functional equation by Dirichlet polynomials having a more convenient, ``factored'' form, as well as slightly reduced lengths. We give a bit more general result, dealing with Dirichlet polynomials $\sum_{n \leq \sqrt{t/2\pi}} \frac{1}{n^{\sigma+it}}$ with $\sigma$ possibly slightly different to 1/2, to assist later with the proof of Theorem 2.

\begin{lem3}\label{polyeditlem}
Let $P \leq \sqrt{T}$ be large, and suppose that $1/\log T < \epsilon < 1/10$, say. Then uniformly for all $\frac{1}{2} - \frac{1}{\log T} \leq \sigma \leq \frac{1}{2} + \frac{1}{\log T}$, say, we have
$$ \int_{T}^{2T} \Biggl|\sum_{n \leq \sqrt{t/2\pi}} \frac{1}{n^{\sigma + it}} - \sum_{\substack{m \leq T^{\epsilon}, \\ m \; \text{is} \; P \; \text{smooth}}} \frac{1}{m^{1/2+it}} \sum_{\substack{n \leq T^{1/2 - 2\epsilon}, \\ n \; \text{is} \; P \; \text{rough}}} \frac{1}{n^{\sigma+it}} \Biggr|^2 dt \ll T \log T ( e^{-\epsilon (\log T)/\log P} + \epsilon ) . $$
\end{lem3}

\begin{proof}[Proof of Lemma \ref{polyeditlem}]
The first minor issue that we must address is the fact that the range of summation $n \leq \sqrt{t/2\pi}$ depends mildly on $t$. To handle this, we just note that if $u \geq 0$ then $\sqrt{(t+u)/2\pi} = \sqrt{1 + u/t} \sqrt{t/2\pi} = \sqrt{t/2\pi} + O(u/\sqrt{t}) = \sqrt{t/2\pi} + O(u/\sqrt{T})$. So if we break up the integral over $T \leq t \leq 2T$ into sub-intervals of length $\asymp \sqrt{T}$, on each such interval the range of summation will be essentially constant.

For example, for all $T \leq t \leq T + \sqrt{T}$ we have $\sum_{n \leq \sqrt{t/2\pi}} \frac{1}{n^{\sigma + it}} = \sum_{n \leq \sqrt{T/2\pi}} \frac{1}{n^{\sigma + it}} + O(T^{-1/4})$. Then by the second part of Lemma \ref{meanvaluelem}, noting that numbers of the form $mn$ with $P$-smooth part $m \leq T^{\epsilon}$ and $P$-rough part $n \leq T^{1/2 - 2\epsilon}$ are a subset of all numbers $n \leq T^{1/2 - \epsilon} \leq \sqrt{T/2\pi}$, we have
\begin{eqnarray}
&& \int_{T}^{T + \sqrt{T}} \Biggl|\sum_{n \leq \sqrt{T/2\pi}} \frac{1}{n^{\sigma + it}} - \sum_{\substack{m \leq T^{\epsilon}, \\ m \; \text{is} \; P \; \text{smooth}}} \frac{1}{m^{1/2+it}} \sum_{\substack{n \leq T^{1/2 - 2\epsilon}, \\ n \; \text{is} \; P \; \text{rough}}} \frac{1}{n^{\sigma+it}} \Biggr|^2 dt \nonumber \\
& \ll & \int_{T}^{T + \sqrt{T}} \Biggl|\sum_{n \leq \sqrt{T/2\pi}} \frac{1}{n^{\sigma + it}} - \sum_{\substack{m \leq T^{\epsilon}, \\ P \; \text{smooth}}} \frac{1}{m^{\sigma+it}} \sum_{\substack{n \leq T^{1/2 - 2\epsilon}, \\ P \; \text{rough}}} \frac{1}{n^{\sigma+it}} \Biggr|^2 dt + \nonumber \\
&& + \int_{T}^{T + \sqrt{T}} \Biggl|\sum_{\substack{m \leq T^{\epsilon}, \\ P \; \text{smooth}}} \frac{1}{m^{\sigma+it}} - \sum_{\substack{m \leq T^{\epsilon}, \\ P \; \text{smooth}}} \frac{1}{m^{1/2+it}} \Biggr|^2 \Biggl| \sum_{\substack{n \leq T^{1/2 - 2\epsilon}, \\ P \; \text{rough}}} \frac{1}{n^{\sigma+it}} \Biggr|^2 dt \nonumber \\
& \ll & \sqrt{T} ( \sum_{n \leq \sqrt{T/2\pi}} \frac{1}{n^{2\sigma}} - \sum_{\substack{m \leq T^{\epsilon}, \\ P \; \text{smooth}}} \frac{1}{m^{2\sigma}} \sum_{\substack{n \leq T^{1/2 - 2\epsilon}, \\ P \; \text{rough}}} \frac{1}{n^{2\sigma}} ) + \sqrt{T} \sum_{\substack{m \leq T^{\epsilon}, \\ P \; \text{smooth}}} \frac{|m^{1/2-\sigma} - 1|^2}{m} \sum_{\substack{n \leq T^{1/2 - 2\epsilon}, \\ P \; \text{rough}}} \frac{1}{n^{2\sigma}}  . \nonumber
\end{eqnarray}
In the first difference of sums here, if a number is $\leq T^{1/2 - 2\epsilon}$ and its $P$-smooth part is $\leq T^{\epsilon}$ then its contribution certainly cancels out, so for this first part we have an upper bound
\begin{eqnarray}
& \ll & \sqrt{T} ( \sum_{T^{1/2 - 2\epsilon} < n \leq \sqrt{T/2\pi}} \frac{1}{n^{2\sigma}} + \sum_{\substack{T^{\epsilon} < m \leq T^{1/2 - 2\epsilon}, \\ m \; \text{is} \; P \; \text{smooth}}} \frac{1}{m^{2\sigma}} \sum_{\substack{n \leq T^{1/2 - 2\epsilon}, \\ n \; \text{is} \; P \; \text{rough}}} \frac{1}{n^{2\sigma}} ) \nonumber \\
& \ll & \sqrt{T} ( \sum_{T^{1/2 - 2\epsilon} < n \leq \sqrt{T/2\pi}} \frac{1}{n} + \sum_{\substack{m > T^{\epsilon}, \\ m \; \text{is} \; P \; \text{smooth}}} \frac{1}{m} \sum_{\substack{n \leq T^{1/2 - 2\epsilon}, \\ n \; \text{is} \; P \; \text{rough}}} \frac{1}{n} ) \ll \sqrt{T} (\epsilon\log T + \frac{\log T}{\log P} \sum_{\substack{m > T^{\epsilon}, \\ m \; \text{is} \; P \; \text{smooth}}} \frac{1}{m}) . \nonumber
\end{eqnarray}
Here we used the estimate $\sum_{\substack{n \leq T^{1/2 - 2\epsilon}, \\ n \; \text{is} \; P \; \text{rough}}} \frac{1}{n} \leq \prod_{P < p \leq T^{1/2 - 2\epsilon}} (1 - \frac{1}{p})^{-1} \ll \frac{\log T}{\log P}$. And we can upper bound $\sum_{\substack{m > T^{\epsilon}, \\ m \; \text{is} \; P \; \text{smooth}}} \frac{1}{m}$ by $T^{-\epsilon/\log P} \sum_{\substack{m : \\ m \; \text{is} \; P \; \text{smooth}}} \frac{1}{m^{1-1/\log P}}$, which is equal to $T^{-\epsilon/\log P} \prod_{p \leq P} (1 - \frac{1}{p^{1-1/\log P}})^{-1} \ll T^{-\epsilon/\log P} \log P$.

In the other sums, since $|1/2 - \sigma| \leq 1/\log T$ we have $|m^{1/2-\sigma} - 1| \ll |1/2 - \sigma|\log m \leq \epsilon$, so the contribution from these is $\ll \epsilon^{2} \sqrt{T} \log T$.

The lemma follows by summing our bounds over all the subintervals of length $\sqrt{T}$ that we broke the integral into.
\end{proof}

\vspace{12pt}
Another important tool for our work will be a construction of smooth functions that approximate the indicator function of an interval fairly well, and can be expanded as a polynomial series with good control on the size of the coefficients (so that we will be able to use Lemmas \ref{meanvaluelem} and \ref{highmomentslem} to investigate the average behaviour of such a function applied to a Dirichlet polynomial).

\begin{approxres1}
For any $R \geq 0$ and $\delta > 0$, there exists a function $\gamma : \R \rightarrow [0,1+\delta]$ with the following properties:
\begin{enumerate}
\item $\gamma(x) \geq 1$ for all $|x| \leq R$;

\item $\gamma(x) \leq \delta$ for all $|x| > R + 1$;

\item for all $l \in \N$ and all $x \in \R$, we have the derivative estimate $|\frac{d^{l}}{dx^{l}} \gamma(x)| \leq \frac{(2R+1) (1+\delta)}{\pi (l + 1)} (\frac{2\pi}{\delta})^{l+1}$.

\end{enumerate}
\end{approxres1}

\begin{proof}[Proof of Approximation Result 1]
Let $b(x)$ be a Beurling--Selberg function majorising the indicator function $\textbf{1}_{|x| \leq 1/2}$, with Fourier transform supported on $[-1/\delta,1/\delta]$. See e.g. Vaaler's paper~\cite{vaaler} for details of the construction and properties of such majorants. Thus we get $b(x) \geq \textbf{1}_{|x| \leq 1/2}$ for all $x \in \R$; and $\int_{-\infty}^{\infty} b(x) dx = 1 + \delta$; and $b(x) = \int_{-1/\delta}^{1/\delta} \hat{b}(t) e^{2\pi i xt} dt$ for all $x \in \R$, where $|\hat{b}(t)| = |\int b(x) e^{-2\pi i x t} dx| \leq 1 + \delta$.

Now we simply define $\gamma(x)$ to be a suitable convolution of $b$, namely
$$ \gamma(x) = \int_{-\infty}^{\infty} \textbf{1}_{|u| \leq R + 1/2} b(x-u) du = \int_{-\infty}^{\infty} \textbf{1}_{|x-u| \leq R + 1/2} b(u) du . $$
Then it is clear that $0 \leq \gamma(x) \leq 1+\delta$ for all $x$, since $b(x)$ is non-negative and its integral over the whole real line is $1+\delta$. Furthermore, we can write 
$$ \gamma(x) = \int_{-\infty}^{\infty} \textbf{1}_{|u| \leq R + 1/2} \textbf{1}_{|x - u| \leq 1/2} du + \int_{-\infty}^{\infty} \textbf{1}_{|u| \leq R + 1/2} (b(x-u) - \textbf{1}_{|x - u| \leq 1/2}) du . $$
When $|x| \leq R$, then for all $u$ satisfying $|x - u| \leq 1/2$ we have $|u| \leq R+1/2$, so the first integral here is equal to 1. And we always have $b(x-u) \geq \textbf{1}_{|x-u| \leq 1/2}$, so the second integral is non-negative and we get that $\gamma(x) \geq 1$, which proves the first statement in the result. When $x \notin [-(R+1),(R+1)]$ the first integral vanishes, since we cannot have $|u| \leq R+1/2$ and $|x-u| \leq 1/2$ if $|x| > R+1$. Thus when $x \notin [-(R+1),(R+1)]$ we have
$$ |\gamma(x)| \leq \int_{-\infty}^{\infty} |b(x-u) - \textbf{1}_{|x - u| \leq 1/2}| du = \delta , $$
since we always have $b(x-u) \geq \textbf{1}_{|x-u| \leq 1/2}$, and $\int_{-\infty}^{\infty} (b(x-u) - \textbf{1}_{|x-u| \leq 1/2}) du = \int_{-\infty}^{\infty} b(x-u) du - 1 = \delta$. So we have proved the second statement in the result as well.

To prove the final statement, we just note that for all $l \in \N$ and $x \in \R$ we have $|\frac{d^{l}}{dx^{l}} b(x)| = |\int_{-1/\delta}^{1/\delta} \hat{b}(t) (2\pi i t)^{l} e^{2\pi i xt} dt| \leq \frac{2^{l+1} \pi^{l} (1+\delta)}{l + 1} (\frac{1}{\delta})^{l+1}$, and therefore
$$ |\frac{d^{l}}{dx^{l}} \gamma(x)| = |\int_{-\infty}^{\infty} \textbf{1}_{|u| \leq R + 1/2} \frac{d^{l}}{dx^{l}} b(x-u) du| \leq \frac{(2R+1) (1+\delta)}{\pi (l + 1)} (\frac{2\pi}{\delta})^{l+1} . $$
\end{proof}

\section{Proof of Theorem 1}\label{thm1proofsec}
As described in the Introduction, the overall structure of this proof will closely model the proof of the corresponding random result, namely the upper bound part of Theorem 1 of Harper~\cite{harperrmflow}. Nevertheless, for the sake of clarity we will give fairly complete details.

Let $P$ be a large quantity, which we will later fix in terms of $T$. For each point $|h| \leq 1/2$, we will wish to approximate $h$ on various different ``scales'' (i.e. to various different levels of precision) so that we can discretise the set of $h$ under consideration at certain points in the argument. Thus we set $h(-1) = h$, and for each $0 \leq j \leq \log\log P - 1$ define
$$ h(j) := \max\{u \leq h(j-1): u = \frac{n}{((\log P)/e^j) \log((\log P)/e^j)} \; \text{for some} \; n \in \Z\} . $$
Corresponding to these different scales, we define the partial Euler products
$$ I_l(h) = I_{l,t}(h) := \prod_{P^{e^{-(l+1)}} < p \leq P^{e^{-l}}} (1 - \frac{1}{p^{1/2+it+ih}})^{-1} , $$
where $t \in \R$. Towards the end of the proof, we will make a comparison with the analogous random products, namely
$$ I_{l,\text{rand}}(h) := \prod_{P^{e^{-(l+1)}} < p \leq P^{e^{-l}}} (1 - \frac{f(p)}{p^{1/2+ih}})^{-1} , $$
where $f(p)$ is a sequence of independent Steinhaus random variables.

Next, we let $g(j) :=  C\min\{\sqrt{\log\log P}, \frac{1}{1-q} \} + 2\log\log(\frac{\log P}{e^j})$, and let $\mathcal{G} = \mathcal{G}_{t}$ denote the ``good'' event that for all $|h| \leq 1/2$ and all $0 \leq j \leq \log\log P - B - 1$, we have
$$ \Biggl( \frac{\log P}{e^j} e^{g(j)} \Biggr)^{-1} \leq \prod_{l = j}^{\lfloor \log\log P \rfloor - B - 1} |I_{l}(h(l))| \leq \frac{\log P}{e^j} e^{g(j)} . $$
Here $B \in \N$ is a certain large constant coming from the probabilistic argument in Proposition 5 of Harper~\cite{harperrmflow} (see the proof of Lemma \ref{mainprob1lem} below), and $C$ is a large constant that we shall fix shortly. The reader who is comparing with \cite{harperrmflow} should note that all this corresponds to the set-up in section 4 there (in slightly modified notation), with $k=0$ and with $x$ replaced by $P^e$.

\vspace{12pt}
We then have the following key estimates.

\begin{keyprop1}
Uniformly for any $P \leq T^{1/(\log\log T)^6}$ that are sufficiently large in terms of $C$, and $\frac{\log P \log\log\log P}{\log T} < \epsilon < 1/10$, and $2/3 \leq q \leq 1$, we have
\begin{eqnarray}
&& \frac{1}{T} \int_{T}^{2T} \textbf{1}_{\mathcal{G}_{t}} \Biggl( \int_{-1/2}^{1/2} |\sum_{\substack{m \leq T^{\epsilon}, \\ m \; \text{is} \; P \; \text{smooth}}} \frac{1}{m^{1/2+i(t+h)}}|^2 |\sum_{\substack{n \leq T^{1/2 - 2\epsilon}, \\ n \; \text{is} \; P \; \text{rough}}} \frac{1}{n^{1/2+i(t+h)}}|^2 dh \Biggr)^q dt \nonumber \\
& \ll & \Biggl( C \log T \min\{1, \frac{1}{(1-q)\sqrt{\log\log P}} \} \Biggr)^q , \nonumber
\end{eqnarray}
where $\textbf{1}$ denotes the indicator function.
\end{keyprop1}

\begin{keyprop2}
Uniformly for any $P \leq T^{1/(\log\log T)^2}$ that are sufficiently large in terms of $C$, and any $2/3 \leq q \leq 1$, we have
$$ \frac{1}{T} \text{meas}\{T \leq t \leq 2T : \mathcal{G}_{t} \; \text{fails}\} \ll e^{-2C\min\{\sqrt{\log\log P}, \frac{1}{1-q}\}} . $$
\end{keyprop2}

These statements are precisely analogous (apart from the conditions on the sizes of $P$ and $\epsilon$, which don't arise there) to the corresponding ones from the random case~\cite{harperrmflow}, although as we shall see some different technical issues arise in the proofs. But using the same recursive procedure as in the random case, we can combine Key Propositions 1 and 2 to prove Theorem 1.

\begin{proof}[Proof of Theorem 1, assuming Key Propositions 1 and 2]
We now fix the choices $\epsilon = \frac{1}{(\log\log T)^2}$ and $P = T^{1/(\log\log T)^8}$, which certainly satisfy the conditions of Key Propositions 1 and 2. (Many other choices of $P, \epsilon$ would also work for this proof, in addition to satisfying the conditions of the propositions we only require that $\epsilon \ll \frac{1}{\sqrt{\log\log T}}$ and $\log\log P \gg \log\log T$, so that the bounds from the propositions and from Lemma \ref{polyeditlem} are ultimately strong enough in terms of $T$. These particular values of $\epsilon, P$ will also work later for the proof of Theorem 2.)

Uniformly for all large $T$ and $0 \leq q \leq 1$, by the approximate functional equation (Zeta Function Result 1) we have
$$ \frac{1}{T} \int_{T}^{2T} \left( \int_{-1/2}^{1/2} |\zeta(1/2 + i(t+h))|^2 dh \right)^q dt \ll \frac{1}{T} \int_{T}^{2T} \Biggl( \int_{-1/2}^{1/2} |\sum_{n \leq \sqrt{\frac{t+h}{2\pi}}} \frac{1}{n^{1/2 + i(t+h)}}|^2 dh \Biggr)^q dt + T^{-q/2} . $$
The second term on the right hand side is certainly acceptable for Theorem 1, whereas the first term is
\begin{eqnarray}
& \ll & \frac{1}{T} \int_{T}^{2T} \Biggl( \int_{-1/2}^{1/2} |\sum_{n \leq \sqrt{\frac{t+h}{2\pi}}} \frac{1}{n^{1/2 + i(t+h)}} - \sum_{\substack{m \leq T^{\epsilon}, \\ m \; \text{is} \; P \; \text{smooth}}} \frac{1}{m^{1/2+i(t+h)}} \sum_{\substack{n \leq T^{1/2 - 2\epsilon}, \\ n \; \text{is} \; P \; \text{rough}}} \frac{1}{n^{1/2+i(t+h)}}|^2 dh \Biggr)^q dt \nonumber \\
&& + \frac{1}{T} \int_{T}^{2T} \Biggl( \int_{-1/2}^{1/2} \Biggl|\sum_{\substack{m \leq T^{\epsilon}, \\ m \; \text{is} \; P \; \text{smooth}}} \frac{1}{m^{1/2+i(t+h)}} \sum_{\substack{n \leq T^{1/2 - 2\epsilon}, \\ n \; \text{is} \; P \; \text{rough}}} \frac{1}{n^{1/2+i(t+h)}} \Biggr|^2 dh \Biggr)^q dt . \nonumber
\end{eqnarray}
We can control the first term here by applying H\"{o}lder's inequality with exponent $1/q$ to the integral over $t$, and then applying Lemma \ref{polyeditlem}, giving an acceptable bound
$$ \ll (\log T ( e^{-\epsilon (\log T)/\log P} + \epsilon ))^q \ll (\frac{\log T}{(\log\log T)^2})^q . $$

\vspace{12pt}
Next, if $1 - \frac{1}{\sqrt{\log\log T}} \leq q \leq 1$ then, simply applying H\"{o}lder's inequality followed by the second part of Lemma \ref{meanvaluelem}, we find the second term in the previous display is
\begin{eqnarray}
& \leq & \Biggl( \frac{1}{T} \int_{T}^{2T} \int_{-1/2}^{1/2} \Biggl|\sum_{\substack{m \leq T^{\epsilon}, \\ m \; \text{is} \; P \; \text{smooth}}} \frac{1}{m^{1/2+i(t+h)}} \sum_{\substack{n \leq T^{1/2 - 2\epsilon}, \\ n \; \text{is} \; P \; \text{rough}}} \frac{1}{n^{1/2+i(t+h)}} \Biggr|^2 dh dt \Biggr)^q \nonumber \\
& \ll & \Biggl( \int_{-1/2}^{1/2} \sum_{\substack{m \leq T^{\epsilon}, \\ m \; \text{is} \; P \; \text{smooth}}} \frac{1}{m} \sum_{\substack{n \leq T^{1/2 - 2\epsilon}, \\ n \; \text{is} \; P \; \text{rough}}} \frac{1}{n} dh \Biggr)^q \ll \log^{q}T , \nonumber
\end{eqnarray}
which is the Theorem 1 bound on this range of $q$. (We have really just applied H\"{o}lder's inequality and a second moment estimate for the zeta function, in a complicated way.) To handle smaller $q$, for each $1/\sqrt{\log\log T} \leq \delta \leq 1/6$ (say) let us define $R(\delta) = R(\delta,T)$ to be
$$ := \sup_{1-2\delta \leq q \leq 1-\delta} \frac{1}{T} \int_{T}^{2T} (\frac{(1-q) \sqrt{\log\log T}}{\log T} \int_{-1/2}^{1/2} |\sum_{\substack{m \leq T^{\epsilon}, \\ P \; \text{smooth}}} \frac{1}{m^{1/2+i(t+h)}} \sum_{\substack{n \leq T^{1/2 - 2\epsilon}, \\ P \; \text{rough}}} \frac{1}{n^{1/2+i(t+h)}}|^2 dh )^q dt . $$
We shall prove a uniform upper bound $R(\delta) \ll 1$, which will imply the Theorem.

Indeed, we may split up the integral over $t$ into the parts where $\mathcal{G}_t$ holds and where it fails. Applying Key Proposition 1 to the first part, we find  its contribution to $R(\delta)$ is $\ll C^{1-\delta} \leq C$. Note that here we need to use the fact that $\sqrt{\log\log P} \asymp \sqrt{\log\log T}$. To handle the contribution from points $t$ for which $\mathcal{G}_t$ fails, for each $1-2\delta \leq q \leq 1-\delta$ let us set $q' = (1+q)/2$, so that $1-\delta \leq q' \leq 1-\delta/2$. Then by H\"{o}lder's inequality with exponents $q'/(q'-q)$ and $q'/q$, we have
\begin{eqnarray}
&& \frac{1}{T} \int_{T}^{2T} \textbf{1}_{\mathcal{G}_{t} \; \text{fails}} (\frac{(1-q) \sqrt{\log\log T}}{\log T} \int_{-1/2}^{1/2} |\sum_{\substack{m \leq T^{\epsilon}, \\ P \; \text{smooth}}} \frac{1}{m^{1/2+i(t+h)}} \sum_{\substack{n \leq T^{1/2 - 2\epsilon}, \\ P \; \text{rough}}} \frac{1}{n^{1/2+i(t+h)}}|^2 dh )^q \nonumber \\
& \leq & (\frac{1}{T} \int_{T}^{2T} \textbf{1}_{\mathcal{G}_{t} \; \text{fails}} dt)^{\frac{q'-q}{q'}} \cdot \nonumber \\
&& \cdot \Biggl(\frac{1}{T} \int_{T}^{2T} (\frac{(1-q) \sqrt{\log\log T}}{\log T} \int_{-1/2}^{1/2} |\sum_{\substack{m \leq T^{\epsilon}, \\ P \; \text{smooth}}} \frac{1}{m^{1/2+i(t+h)}} \sum_{\substack{n \leq T^{1/2 - 2\epsilon}, \\ P \; \text{rough}}} \frac{1}{n^{1/2+i(t+h)}}|^2 dh )^{q'} \Biggr)^{q/q'} . \nonumber
\end{eqnarray}
The point is that we have $(q'-q)/q' = (1-q)/(2q') \geq \delta/2$, whilst by Key Proposition 2 we have $\frac{1}{T} \int_{T}^{2T} \textbf{1}_{\mathcal{G}_{t} \; \text{fails}} dt \ll e^{-2C\min\{\sqrt{\log\log P}, \frac{1}{1-q}\}} \leq e^{-C/\delta}$. Meanwhile, in the second bracket we have $1 - q = 2(1-q')$, and we always have $q/q' \leq 1$, so the contribution from this bracket is $\ll R(\delta/2)^{q/q'} \leq (1 + R(\delta/2))$. Putting everything together, we have shown that
$$ R(\delta) \ll C + e^{-C/2} (1 + R(\delta/2)) \ll C + e^{-C/2}R(\delta/2) . $$

Iterating this recursive bound (with $C$ {\em fixed} sufficiently large to compensate for the implicit constant), replacing $\delta$ by $\delta/2, \delta/4, \delta/8$, etc., we see that uniformly for $1/\sqrt{\log\log T} \leq \delta \leq 1/6$ we have
$$ R(\delta) \ll 1 + R(1/\sqrt{\log\log T}) . $$
Our earlier calculations for the range $1 - \frac{1}{\sqrt{\log\log T}} \leq q \leq 1$ imply that $R(1/\sqrt{\log\log T}) \ll 1$, so we have our uniform upper bound $R(\delta) \ll 1$, as required.
\end{proof}

\vspace{12pt}
It remains to prove our Key Propositions. We begin with the easier of the two, namely Key Proposition 2. The proof is ultimately just an application of the union bound (exploiting the discrete nature of the set of approximating points $h(l)$), as in the random case~\cite{harperrmflow}. However, whereas in the random case one obtains the tail estimates in the argument by computing the second moment of random Euler products, here this would be more difficult to achieve (one would first need to approximate the products by sums of suitable length) and instead we shall apply Lemma \ref{highmomentslem} directly to prime number sums.  

\begin{proof}[Proof of Key Proposition 2]
Using the definition of $\mathcal{G}_t$ and the union bound, we can upper bound the left hand side in the proposition by
$$ \sum_{j=0}^{\lfloor \log\log P \rfloor - B - 1} \sum_{h(j)} \frac{1}{T} \text{meas}\Biggl\{T \leq t \leq 2T : \Biggl|\sum_{l=j}^{\lfloor \log\log P \rfloor - B - 1} \log|I_{l,t}(h(l))| \Biggr| > \log\log P - j + g(j) \Biggr\} , $$
where the second sum is over all possible values of $h(j)$ as $|h| \leq 1/2$ varies. (Note there are $\asymp ((\log P)/e^j) \log((\log P)/e^j)$ such values, and once we know $h(j)$ this uniquely determines the values $h(l)$ for all $l \geq j$.) And by definition of $I_{l,t}(h)$, the sum over $l$ is
\begin{eqnarray}
& = & - \Re \sum_{l=j}^{\lfloor \log\log P \rfloor - B - 1} \sum_{P^{e^{-(l+1)}} < p \leq P^{e^{-l}}} \log(1 - \frac{1}{p^{1/2+it+ih(l)}}) \nonumber \\
& = & \Re \sum_{l=j}^{\lfloor \log\log P \rfloor - B - 1} \sum_{P^{e^{-(l+1)}} < p \leq P^{e^{-l}}} (\frac{1}{p^{1/2+it+ih(l)}} + \frac{1}{2p^{1+2it+2ih(l)}} + O(\frac{1}{p^{3/2}})) . \nonumber
\end{eqnarray}
So in order for this to have absolute value $> \log\log P - j + g(j)$, we must have
$$ \Biggl|\sum_{l=j}^{\lfloor \log\log P \rfloor - B - 1} \sum_{\substack{P^{e^{-(l+1)}} < p \leq P^{e^{-l}}, \\ v=1,2}} \frac{1}{vp^{v(1/2+it+ih(l))}} \Biggr| > \log\log P - j + g(j) + O(1) . $$

Now for any $k \in \N$ such that $P^{2k e^{-j}} < T$ (so the second part of Lemma \ref{highmomentslem} is applicable), the measure of the set of $T \leq t \leq 2T$ for which this inequality holds may be bounded by
\begin{eqnarray}
&& \frac{1}{(\log\log P - j + g(j) + O(1))^{2k}}  \int_{T}^{2T} \Biggl|\sum_{l=j}^{\lfloor \log\log P \rfloor - B - 1} \sum_{\substack{P^{e^{-(l+1)}} < p \leq P^{e^{-l}}, \\ v=1,2}} \frac{1}{vp^{v(1/2+it+ih(l))}} \Biggr|^{2k} dt \nonumber \\
& \ll & \frac{1}{(\log\log P - j + g(j) + O(1))^{2k}} T (k!) (\sum_{p \leq P^{e^{-j}}} \frac{1}{p} + O(1))^{k} . \nonumber
\end{eqnarray}
Using Stirling's formula together with the Mertens estimate for $\sum 1/p$, this is all
$$ \ll T \sqrt{k} \left(\frac{k(\log\log P - j + O(1))}{e(\log\log P - j + g(j) + O(1))^{2}} \right)^{k} . $$
A quick calculation shows that choosing $k = \lfloor (\log\log P - j + g(j) + O(1))^{2}/(\log\log P - j + O(1)) \rfloor$ is roughly optimal here, and yields an overall bound
$$ \\ T \sqrt{k} \exp\left\{- \lfloor \frac{(\log\log P - j + g(j) + O(1))^{2}}{\log\log P - j + O(1)} \rfloor \right\} \ll T \sqrt{k} \frac{e^j}{\log P} e^{-2g(j) - \frac{(g(j) + O(1))^2}{\log\log P - j + O(1)}} . $$
Recalling the definition of $g(j) = C\min\{\sqrt{\log\log P}, \frac{1}{1-q} \} + 2\log\log(\frac{\log P}{e^j})$, when we have $C\min\{\sqrt{\log\log P}, \frac{1}{1-q} \} \leq (\log\log P - j)$ then we have $k \asymp (\log\log P - j)$, so our bound is $\ll T \frac{e^j}{\log P} \frac{1}{\log^{7/2}((\log P)/e^j)} e^{-2C\min\{\sqrt{\log\log P}, \frac{1}{1-q} \}}$. And in the other case where $C\min\{\sqrt{\log\log P}, \frac{1}{1-q} \} > (\log\log P - j)$, we have at least as good a bound thanks to the additional term $- \frac{(g(j) + O(1))^2}{\log\log P - j + O(1)}$ in the exponential. Notice also that the condition $P^{2k e^{-j}} < T$ will always be satisfied because of our assumption that $P \leq T^{1/(\log\log T)^2}$.

Putting everything together, we have shown that
$$ \frac{1}{T} \text{meas}\{T \leq t \leq 2T : \mathcal{G}_{t} \; \text{fails}\} \ll \sum_{j=0}^{\lfloor \log\log P \rfloor - B - 1} \sum_{h(j)} \frac{e^j}{\log P} \frac{1}{\log^{7/2}((\log P)/e^j)} e^{-2C\min\{\sqrt{\log\log P}, \frac{1}{1-q} \}} . $$
For each $j$, the number of values $h(j)$ that we sum over is $\ll \frac{\log P}{e^j}  \log((\log P)/e^j)$. Inserting this bound proves the proposition.
\end{proof}

Now we shall prove Key Proposition 1, which is where most of our technical work with approximating functions and polynomial expansions will arise. Towards the end of the proof we will manage to pass to the random case, and will need the following lemma which imports a key probabilistic estimate of Harper~\cite{harperrmflow}. For each $|h| \leq 1/2$, let $\mathcal{G}'(h)$ denote the event that for all $0 \leq j \leq \log\log P - B - 1$, we have
$$ \Biggl( \frac{\log P}{e^j} e^{g(j) + 19/3} \Biggr)^{-1} \leq \prod_{l = j}^{\lfloor \log\log P \rfloor - B - 1} |I_{l, \text{rand}}(h(l))| \leq \frac{\log P}{e^j} e^{g(j) + 19/3} , $$
where $I_{l, \text{rand}}(\cdot)$ and $h(l)$ are the random Euler products and approximating points that we introduced earlier. (The insertion of the additional factor 19/3 in the exponents is to absorb some shifts by constants when approximating the Euler products, as will soon become clear in the proof of Key Proposition 1.)  

\begin{lem4}\label{mainprob1lem}
Let $f(n)$ denote a Steinhaus random multiplicative function. Uniformly for all large $P \leq \sqrt{T}$ and $0 < \epsilon < 1/10$ and $2/3 \leq q \leq 1$ and $|h| \leq 1/2$, we have
\begin{eqnarray}
&& \E \textbf{1}_{\mathcal{G}^{'}(h)} |\sum_{\substack{m \leq T^{\epsilon}, \\ m \; \text{is} \; P \; \text{smooth}}} \frac{f(m)}{m^{1/2+ih}}|^2 |\sum_{\substack{n \leq T^{1/2 - 2\epsilon}, \\ n \; \text{is} \; P \; \text{rough}}} \frac{f(n)}{n^{1/2+ih}}|^2 \nonumber \\
& \ll & \log T \left( C \min\{1, \frac{1}{(1-q)\sqrt{\log\log P}} \} + e^{-\epsilon (\log T)/\log P} \right) . \nonumber
\end{eqnarray}
\end{lem4}

\begin{proof}[Proof of Lemma \ref{mainprob1lem}]
Since the underlying $f(p)$ are independent for distinct primes $p$, and the event $\mathcal{G}^{'}(h)$ only depends on the random variables $(f(p))_{p \leq P}$, our expectation factors as
$$ \E \textbf{1}_{\mathcal{G}^{'}(h)} |\sum_{\substack{m \leq T^{\epsilon}, \\ m \; \text{is} \; P \; \text{smooth}}} \frac{f(m)}{m^{1/2+ih}}|^2 \cdot \E |\sum_{\substack{n \leq T^{1/2 - 2\epsilon}, \\ n \; \text{is} \; P \; \text{rough}}} \frac{f(n)}{n^{1/2+ih}}|^2 \asymp \frac{\log T}{\log P} \E \textbf{1}_{\mathcal{G}^{'}(h)} |\sum_{\substack{m \leq T^{\epsilon}, \\ m \; \text{is} \; P \; \text{smooth}}} \frac{f(m)}{m^{1/2+ih}}|^2 . $$
We also have that $\E \textbf{1}_{\mathcal{G}^{'}(h)} |\sum_{\substack{m \leq T^{\epsilon}, \\ m \; \text{is} \; P \; \text{smooth}}} \frac{f(m)}{m^{1/2+ih}}|^2$ is
\begin{eqnarray}
& \ll & \E \textbf{1}_{\mathcal{G}^{'}(h)} |\prod_{p \leq P} (1 - \frac{f(p)}{p^{1/2+ih}})|^{-2} + \E |\sum_{\substack{m \leq T^{\epsilon}, \\ m \; \text{is} \; P \; \text{smooth}}} \frac{f(m)}{m^{1/2+ih}} - \prod_{p \leq P} (1 - \frac{f(p)}{p^{1/2+ih}})^{-1}|^2 \nonumber \\
& = & \E \textbf{1}_{\mathcal{G}^{'}(h)} |\prod_{p \leq P} (1 - \frac{f(p)}{p^{1/2+ih}})|^{-2} + \sum_{\substack{m > T^{\epsilon}, \\ m \; \text{is} \; P \; \text{smooth}}} \frac{1}{m} , \nonumber
\end{eqnarray}
using the fact that $\prod_{p \leq P} (1 - \frac{f(p)}{p^{1/2+ih}})^{-1} = \sum_{\substack{m = 1, \\ m \; \text{is} \; P \; \text{smooth}}}^{\infty} \frac{f(m)}{m^{1/2+ih}}$ together with the orthogonality of the random variables $f(m)$. As in the proof of Lemma \ref{polyeditlem}, we have $\sum_{\substack{m > T^{\epsilon}, \\ m \; \text{is} \; P \; \text{smooth}}} \frac{1}{m} \ll T^{-\epsilon/\log P} \log P$, which gives an acceptable contribution. And finally applying Proposition 5 from Harper~\cite{harperrmflow} (in the same manner as in the proof of Key Proposition 1 there), we have
\begin{eqnarray}
\E \textbf{1}_{\mathcal{G}^{'}(h)} |\prod_{p \leq P} (1 - \frac{f(p)}{p^{1/2+ih}})|^{-2} & \ll & C \min\{1, \frac{1}{(1-q)\sqrt{\log\log P}} \} \E |\prod_{p \leq P} (1 - \frac{f(p)}{p^{1/2+ih}})|^{-2} \nonumber \\
& \ll & C \min\{1, \frac{1}{(1-q)\sqrt{\log\log P}} \} \log P , \nonumber
\end{eqnarray}
using that $\E |\prod_{p \leq P} (1 - \frac{f(p)}{p^{1/2+ih}})|^{-2} = \sum_{\substack{m = 1, \\ m \; \text{is} \; P \; \text{smooth}}}^{\infty} \frac{1}{m} \asymp \log P$.
\end{proof}

\begin{proof}[Proof of Key Proposition 1]
Applying H\"{o}lder's Inequality with exponent $1/q$ to the integral over $t$, we find it will suffice to prove that uniformly for all $|h| \leq 1/2$, we have
$$ \frac{1}{T} \int_{T}^{2T} \textbf{1}_{\mathcal{G}_{t}} |\sum_{\substack{m \leq T^{\epsilon}, \\ P \; \text{smooth}}} \frac{1}{m^{1/2+i(t+h)}}|^2 |\sum_{\substack{n \leq T^{1/2 - 2\epsilon}, \\ P \; \text{rough}}} \frac{1}{n^{1/2+i(t+h)}}|^2 dt \ll C \log T \min\{1, \frac{1}{(1-q)\sqrt{\log\log P}} \} . $$

We can upper bound $\textbf{1}_{\mathcal{G}_{t}}$ here by $\textbf{1}_{\mathcal{G}_{t}(h)}$, where $\mathcal{G}_{t}(h)$ denotes the event that for the specific $h$ we are looking at (rather than all $|h| \leq 1/2$) we have
$$ \Biggl( \frac{\log P}{e^j} e^{g(j)} \Biggr)^{-1} \leq \prod_{l = j}^{\lfloor \log\log P \rfloor - B - 1} |I_{l}(h(l))| \leq \frac{\log P}{e^j} e^{g(j)}  $$
for all $0 \leq j \leq \log\log P - B - 1$. Taking logarithms as usual, and recalling the definition of $I_{l}(h)$, we have
\begin{eqnarray}
\log \prod_{l = j}^{\lfloor \log\log P \rfloor - B - 1} |I_{l}(h(l))| & = & \Re \sum_{l = j}^{\lfloor \log\log P \rfloor - B - 1} \log(I_{l}(h(l))) \nonumber \\
& = & \Re \sum_{l = j}^{\lfloor \log\log P \rfloor - B - 1} \sum_{P^{e^{-(l+1)}} < p \leq P^{e^{-l}}} \sum_{v=1}^{\infty} \frac{1}{v} \frac{1}{p^{v(1/2+it+ih(l))}} . \nonumber
\end{eqnarray}
The contribution from all the terms here with $v \geq 3$ is trivially at most $\sum_{p} \sum_{v \geq 3} \frac{1}{vp^{v/2}} \leq \frac{4}{3} \sum_{p} \frac{1}{p^{3/2}} \leq \frac{4}{3} \int_{1}^{\infty} \frac{dw}{w^{3/2}} = 8/3$. So letting $\gamma_{j}(\cdot)$ denote the function supplied by Approximation Result 1, with the choice
$$ R = R_j = \log\log P - j + g(j) + 8/3 $$
(and with the choice of small $\delta > 0$ to be determined at the end of the proof), we can further upper bound $\textbf{1}_{\mathcal{G}_{t}(h)}$ by
\begin{equation}\label{gammajprod}
\prod_{j=0}^{\lfloor \log\log P \rfloor - B - 1} \gamma_{j}\Biggl(\Re \sum_{l=j}^{\lfloor \log\log P \rfloor - B - 1} \sum_{\substack{P^{e^{-(l+1)}} < p \leq P^{e^{-l}}, \\ v = 1,2}} \frac{1}{v p^{v(1/2+it+ih(l))}} \Biggr) .
\end{equation}
Note also that by part (iii) of Approximation Result 1, for any $k \in \N$ (to be fixed later) we can write $\gamma_{j}(\cdot) = \gamma^{\text{main}}_{j}(\cdot) + \gamma^{\text{error}}_{j}(\cdot)$, where $\gamma^{\text{main}}_{j}(\cdot)$ denotes the degree $2k-1$ Taylor expansion of $\gamma_{j}(\cdot)$ about 0, and $\gamma^{\text{error}}_{j}(\cdot)$ satisfies $|\gamma^{\text{error}}_{j}(\cdot)| \leq \frac{(2R_j + 1)(1 + \delta)}{2\pi k} (\frac{2\pi}{\delta})^{2k+1} \frac{|\cdot|^{2k}}{(2k)!} \leq \frac{3 \log\log P}{k \delta} (\frac{e \pi}{k \delta})^{2k} |\cdot|^{2k}$. (Here we used the bound $(2k!) \geq (2k/e)^{2k}$.)

Now expanding out the product over $j$ in \eqref{gammajprod}, using the fact that $|\gamma_{j}(\cdot)| \leq 1 + \delta$ and so $|\gamma^{\text{main}}_{j}(\cdot)| \leq 1 + \delta + |\gamma^{\text{error}}_{j}(\cdot)|$, the contribution from everything involving an error term $\gamma^{\text{error}}_{j}(\cdot)$ is crudely at most
\begin{eqnarray}
&& \sum_{j=0}^{\lfloor \log\log P \rfloor - B - 1} |\gamma^{\text{error}}_{j}(\Re \sum_{l} \sum_{p, v} \frac{1}{vp^{v(1/2+it+ih(l))}})| \prod_{k \neq j}\Biggl(1 + \delta + |\gamma^{\text{error}}_{k}(\Re \sum_{l} \sum_{p, v} \frac{1}{vp^{v(1/2+it+ih(l))}} )| \Biggr) \nonumber \\
& \leq & \log\log P \cdot \max_{j} |\gamma^{\text{error}}_{j}(\Re \sum_{l,p,v} \frac{1}{vp^{v(1/2+it+ih(l))}})| \Biggl(1 + \delta + |\gamma^{\text{error}}_{j}(\Re \sum_{l,p,v} \frac{1}{vp^{v(1/2+it+ih(l))}} )| \Biggr)^{\lfloor \log\log P \rfloor} . \nonumber
\end{eqnarray}
Using our upper bound for $|\gamma^{\text{error}}_{j}(\cdot)|$, we can further bound the contribution from all these error terms by
\begin{eqnarray}
& \ll & \frac{(\log\log P)^2}{k \delta} \max_{j} (\frac{e \pi}{k \delta})^{2k} |\sum_{l,p,v} \frac{1}{vp^{v(1/2+it+ih(l))}}|^{2k} \cdot \nonumber \\
&& \cdot \Biggl(1 + \delta + \frac{3 \log\log P}{k \delta} (\frac{e \pi}{k \delta})^{2k} |\sum_{l,p,v} \frac{1}{vp^{v(1/2+it+ih(l))}}|^{2k} \Biggr)^{\lfloor \log\log P \rfloor} . \nonumber
\end{eqnarray}
Recall here that the range of summation over $l$ depends on $j$, so the maximum over $j$ cannot (yet) be dispensed with. Provided we have $k \geq \frac{100\log\log P}{\delta}$, say, we can tidy things up a bit and obtain the bound
\begin{eqnarray}
& \ll & (\log\log P) \max_{j} (\frac{e \pi}{k \delta})^{2k} |\sum_{l,p,v} \frac{1}{vp^{v(1/2+it+ih(l))}}|^{2k} \Biggl(1 + \delta + (\frac{e \pi}{k \delta})^{2k} |\sum_{l,p,v} \frac{1}{vp^{v(1/2+it+ih(l))}}|^{2k} \Biggr)^{\lfloor \log\log P \rfloor} \nonumber \\
& \leq & (\log\log P) \sum_{j=0}^{\lfloor \log\log P \rfloor - B - 1} (\frac{e \pi}{k \delta})^{2k} |\sum_{l,p,v} \frac{1}{vp^{v(1/2+it+ih(l))}}|^{2k} 2^{\lfloor \log\log P \rfloor} \Biggl( (1 + \delta)^{\lfloor \log\log P \rfloor} + \nonumber \\
&& + (\frac{e \pi}{k \delta})^{2k\lfloor \log\log P \rfloor} |\sum_{l,p,v} \frac{1}{vp^{v(1/2+it+ih(l))}}|^{2k\lfloor \log\log P \rfloor} \Biggr) . \nonumber
\end{eqnarray}

Inserting all this in the original $t$-integral, we can use the first part of Lemma \ref{highmomentslem} (which is applicable provided we ultimately have $T^{1/2 - \epsilon} P^{2k(\lfloor \log\log P \rfloor + 1)} < T$) to bound the integrated contribution from all the error terms. Noting that if $m$ is $P$-smooth and $n$ is $P$-rough and $\mathcal{P}$ is a subset of the primes $\leq P$, then the function $\tilde{d}(\cdot)$ from Lemma \ref{highmomentslem} satisfies $\tilde{d}(mn) = \tilde{d}(m) \leq d(m)$ (the classical divisor function), we obtain overall that uniformly for any $|h| \leq 1/2$, we have
\begin{eqnarray}\label{bigapproxdisplay}
&& \frac{1}{T} \int_{T}^{2T} \textbf{1}_{\mathcal{G}_{t}} |\sum_{\substack{m \leq T^{\epsilon}, \\ m \; \text{is} \; P \; \text{smooth}}} \frac{1}{m^{1/2+i(t+h)}} \sum_{\substack{n \leq T^{1/2 - 2\epsilon}, \\ n \; \text{is} \; P \; \text{rough}}} \frac{1}{n^{1/2+i(t+h)}}|^2 dt \nonumber \\
& \ll & \frac{1}{T} \int_{T}^{2T} \prod_{j=0}^{\lfloor \log\log P \rfloor - B - 1} \gamma_{j}^{\text{main}}(\Re \sum_{l,p,v} \frac{1}{v p^{v(1/2+it+ih(l))}} ) |\sum_{\substack{m \leq T^{\epsilon}, \\  P \; \text{smooth}}} \frac{1}{m^{1/2+i(t+h)}} \sum_{\substack{n \leq T^{1/2 - 2\epsilon}, \\  P \; \text{rough}}} \frac{1}{n^{1/2+i(t+h)}}|^2 dt \nonumber \\
&& + (\log\log P)^2 2^{\log\log P} (\frac{e\pi}{k\delta})^{2k} \Biggl( \sum_{\substack{m \leq T^{\epsilon}, \\ P \; \text{smooth}}} \frac{d(m)}{m} \sum_{\substack{n \leq T^{1/2 - 2\epsilon}, \\ P \; \text{rough}}} \frac{1}{n} \Biggr) \Biggl((1+\delta)^{\log\log P} (k!) (\sum_{p \leq P} \frac{2}{p} + O(1))^k \nonumber \\
&& + (\frac{e\pi}{k\delta})^{2k\lfloor \log\log P \rfloor} (k(\lfloor \log\log P \rfloor + 1))! (\sum_{p \leq P} \frac{2}{p} + O(1))^{k(\lfloor \log\log P \rfloor + 1)} \Biggr) .
\end{eqnarray}
To estimate the scary looking second term in \eqref{bigapproxdisplay}, note that $\sum_{\substack{m \leq T^{\epsilon}, \\ P \; \text{smooth}}} \frac{d(m)}{m} \sum_{\substack{n \leq T^{1/2 - 2\epsilon}, \\ P \; \text{rough}}} \frac{1}{n} \ll \prod_{p \leq P} (1 - \frac{1}{p})^{-2} \cdot \frac{\log T}{\log P} \ll \log P \cdot \log T$. Combining this with Stirling's formula and the Mertens estimate, and the upper bound $(1+\delta)^{\log\log P} \leq \log^{\delta}P$, we find these second sums are all
\begin{eqnarray}
& \ll & (\log\log P)^2 2^{\log\log P} \log P \log T \Biggl((\log^{\delta}P) \sqrt{k} \Biggl(\pi \sqrt{\frac{e}{k}} \frac{\sqrt{2(\log\log P + O(1))}}{\delta} \Biggr)^{2k} + \nonumber \\
&& + \sqrt{k\log\log P} \Biggl(\pi \sqrt{\frac{e}{k}} \frac{\sqrt{2} (\log\log P + O(1))}{\delta} \Biggr)^{2k(\lfloor \log\log P \rfloor + 1)} \Biggr) . \nonumber
\end{eqnarray}
In particular, if we take $k = \lfloor (\frac{100\log\log P}{\delta})^2 \rfloor$, say, this is all
\begin{equation}\label{mainerrorkp1}
\ll \log T (\log\log P)^{5/2} (\log^{1 + \log 2 + \delta}P) \sqrt{k} e^{-2k} ,
\end{equation}
which is more than good enough for the proposition.

Finally, since $\gamma_{j}^{\text{main}}(\cdot)$ is just a polynomial of degree $2k-1$, whose degree $d$ coefficient has absolute value at most $\frac{(2R_j + 1)(1 + \delta)}{\pi (d+1)} (\frac{2\pi}{\delta})^{d+1} \frac{1}{d!} \leq (\log\log P) (\frac{2 \pi}{\delta})^{d+1} \frac{1}{(d+1)!}$, we can apply the first part of Lemma \ref{meanvaluelem} to calculate the remaining integral in \eqref{bigapproxdisplay}. We find this is equal to
$$ \E \prod_{j=0}^{\lfloor \log\log P \rfloor - B - 1} \gamma_{j}^{\text{main}}(\Re \sum_{l,p,v} \frac{f(p)^v}{v p^{v(1/2+ih(l))}} ) |\sum_{\substack{m \leq T^{\epsilon}, \\ m \; \text{is} \; P \; \text{smooth}}} \frac{f(m)}{m^{1/2+ih}}|^2 |\sum_{\substack{n \leq T^{1/2 - 2\epsilon}, \\ n \; \text{is} \; P \; \text{rough}}} \frac{f(n)}{n^{1/2+ih}}|^2 , $$
where $f(n)$ is a Steinhaus random multiplicative function, up to a remainder term that has order at most
$$ \frac{P^{8k} T^{1/2 - \epsilon}}{T} \sum_{\substack{m \leq T^{\epsilon}, \\ P \; \text{smooth}}} \frac{d(m)}{m} \sum_{\substack{n \leq T^{1/2 - 2\epsilon}, \\ P \; \text{rough}}} \frac{1}{n} (2k\lfloor \log\log P \rfloor)! \Biggl( \sum_{d=0}^{2k-1} \frac{\log\log P}{(d+1)!} (\frac{2 \pi}{\delta})^{d+1} (\sum_{p \leq P} \frac{2}{p} + O(1))^d \Biggr)^{\log\log P} . $$
The remainder term may be put into this form using the same calculations as in the proof of Lemma \ref{highmomentslem} and as above, on also noting that the largest integers appearing when we expand out the Dirichlet polynomials in $\prod_{j=0}^{\lfloor \log\log P \rfloor - B - 1} \gamma_{j}^{\text{main}}(\Re \sum_{l,p,v} \frac{1}{v p^{v(1/2+it+ih(l))}} )$ have size $\leq \prod_{j=0}^{\lfloor \log\log P \rfloor - B - 1} P^{2(2k-1)e^{-j}} \leq P^{8k}$. Now the Mertens estimate implies that $(\sum_{p \leq P} \frac{2}{p} + O(1))^d = (2\log\log P + O(1))^d$, so using the series expansion of the exponential this remainder is
\begin{eqnarray}\label{othererrorkp1}
& \leq & \frac{P^{8k} T^{1/2 - \epsilon}}{T} \sum_{\substack{m \leq T^{\epsilon}, \\ P \; \text{smooth}}} \frac{d(m)}{m} \sum_{\substack{n \leq T^{1/2 - 2\epsilon}, \\ P \; \text{rough}}} \frac{1}{n} (2k\lfloor \log\log P \rfloor)! \Biggl( \exp\{\frac{4\pi (\log\log P + O(1))}{\delta}\} \Biggr)^{\log\log P} \nonumber \\
& \ll & \frac{P^{8k}}{T^{1/2 + \epsilon}} \log T \log P (2k\lfloor \log\log P \rfloor)! \exp\{\frac{4\pi (\log\log P + O(1))^2}{\delta}\} .
\end{eqnarray}
Meanwhile, by reversing the above calculations we find that up to the same acceptable error term \eqref{mainerrorkp1} obtained earlier, the expectation involving $f(n)$ is equal to
$$ \E \prod_{j=0}^{\lfloor \log\log P \rfloor - B - 1} \gamma_{j}(\Re \sum_{l,p,v} \frac{f(p)^v}{v p^{v(1/2+ih(l))}} ) |\sum_{\substack{m \leq T^{\epsilon}, \\ m \; \text{is} \; P \; \text{smooth}}} \frac{f(m)}{m^{1/2+ih}}|^2 |\sum_{\substack{n \leq T^{1/2 - 2\epsilon}, \\ n \; \text{is} \; P \; \text{rough}}} \frac{f(n)}{n^{1/2+ih}}|^2 , $$
in other words we can replace $\gamma_{j}^{\text{main}}(\cdot)$ by $\gamma_{j}(\cdot)$ again. 

But by construction of the functions $\gamma_j$, we always have $0 \leq \gamma_j(\cdot) \leq 1 + \delta$ and we have $|\gamma_j(\cdot)| \leq \delta$ when $|\cdot| > R_j + 1 = \log\log P - j + g(j) + 11/3$. So we have the upper bound
$$ \prod_{j=0}^{\lfloor \log\log P \rfloor - B - 1} \gamma_{j}( \Re \sum_{l,p,v} \frac{f(p)^v}{v p^{v(1/2+ih(l))}} ) \leq (1 + \delta)^{\log\log P} (\textbf{1}_{\mathcal{G}^{'}}(h) + \delta) , $$
where $\mathcal{G}^{'}(h)$ is (as defined previously) the same as the event $\mathcal{G}_{t}(h)$, except with $g(j)$ replaced by $g(j) + 19/3$ and with $p^{-it}$ replaced by $f(p)$. Here $19/3 = 11/3 + 8/3$ takes account of our trivial estimate 8/3 for the contribution of prime cubes and higher powers. If we finally make the choice $\delta = 1/\log\log P$, then $(1 + \delta)^{\log\log P} \ll 1$ and so our overall expectation is
$$ \ll \E \left(\textbf{1}_{\mathcal{G}^{'}(h)} + \frac{1}{\log\log P} \right) |\sum_{\substack{m \leq T^{\epsilon}, \\ m \; \text{is} \; P \; \text{smooth}}} \frac{f(m)}{m^{1/2+ih}}|^2 |\sum_{\substack{n \leq T^{1/2 - 2\epsilon}, \\ n \; \text{is} \; P \; \text{rough}}} \frac{f(n)}{n^{1/2+ih}}|^2 . $$
Using Lemma \ref{mainprob1lem} and the mean square estimate $\E |\sum_{\substack{m \leq T^{\epsilon}, \\ P \; \text{smooth}}} \frac{f(m)}{m^{1/2+ih}}|^2 |\sum_{\substack{n \leq T^{1/2 - 2\epsilon}, \\ P \; \text{rough}}} \frac{f(n)}{n^{1/2+ih}}|^2 = \sum_{\substack{m \leq T^{\epsilon}, \\ P \; \text{smooth}}} \frac{1}{m} \sum_{\substack{n \leq T^{1/2 - 2\epsilon}, \\ P \; \text{rough}}} \frac{1}{n} \ll \log T$, and recalling that $\epsilon > \frac{\log P \log\log\log P}{\log T}$ in Key Proposition 1, we deduce this is all
$$ \ll \log T \left( C \min\{1, \frac{1}{(1-q)\sqrt{\log\log P}} \} + e^{-\epsilon (\log T)/\log P} \right) \ll C \log T \min\{1, \frac{1}{(1-q)\sqrt{\log\log P}} \} , $$
which is acceptable for the proposition.

It only remains to note that our choice $k = \lfloor (\frac{100\log\log P}{\delta})^2 \rfloor \asymp (\log\log P)^4$ does satisfy our required bound $T^{1/2 - \epsilon} P^{2k(\lfloor \log\log P \rfloor + 1)} < T$ from earlier, and the outstanding remainder term \eqref{othererrorkp1} is by far small enough, in view of our assumption that $P \leq T^{1/(\log\log T)^6}$.
\end{proof}

\section{Further tools}
In addition to the tools we have already introduced, we shall need a couple more to complete the proof of Theorem 2.

Firstly, in Theorem 2 we want to study large values of $\zeta(1/2+it)$, but most of our arguments operate on the level of the Dirichlet polynomial $\sum_{\substack{m \leq T^{\epsilon}, \\ m \; \text{is} \; P \; \text{smooth}}} \frac{1}{m^{1/2+it}}$. So we want to know that if $\zeta(1/2+it)$ is large, then usually $\sum_{\substack{m \leq T^{\epsilon}, \\ m \; \text{is} \; P \; \text{smooth}}} \frac{1}{m^{1/2+it}}$ must be remarkably large as well. We will deduce this from the following fourth moment result.

\begin{lem5}\label{fourthmomentlem}
Uniformly for all large $P \leq \sqrt{T}$ and $0 < \epsilon < 1/10$, and $V > 0$, and $\frac{1}{2} - \frac{1}{\log T} \leq \sigma \leq \frac{1}{2} + \frac{1}{\log T}$, we have
$$ \int_{\substack{T \leq t \leq 2T: \\ |\sum_{\substack{m \leq T^{\epsilon}, \\ P \; \text{smooth}}} \frac{1}{m^{1/2+it}}| \leq (\log P)/V}} \Biggl|\sum_{\substack{m \leq T^{\epsilon}, \\ m \; \text{is} \; P \; \text{smooth}}} \frac{1}{m^{1/2+it}} \sum_{\substack{n \leq T^{1/2 - 2\epsilon}, \\ n \; \text{is} \; P \; \text{rough}}} \frac{1}{n^{\sigma+it}} \Biggr|^4 dt \ll \frac{T}{V^2} \log^{3}T \frac{\log T}{\log P} . $$
\end{lem5}

\begin{proof}[Proof of Lemma \ref{fourthmomentlem}]
We simply upper bound the left hand side by
\begin{eqnarray}
&& \left(\frac{\log P}{V} \right)^2 \int_{T}^{2T} \Biggl|\sum_{\substack{m \leq T^{\epsilon}, \\ m \; \text{is} \; P \; \text{smooth}}} \frac{1}{m^{1/2+it}} \Biggr|^2 \Biggl|\sum_{\substack{n \leq T^{1/2 - 2\epsilon}, \\ n \; \text{is} \; P \; \text{rough}}} \frac{1}{n^{\sigma+it}} \Biggr|^4 dt \nonumber \\
& = & \left(\frac{\log P}{V} \right)^2 \int_{T}^{2T} \Biggl|\sum_{\substack{m \leq T^{\epsilon}, \\ m \; \text{is} \; P \; \text{smooth}}} \frac{1}{m^{1/2+it}} \sum_{\substack{N \leq T^{1 - 4\epsilon}, \\ N \; \text{is} \; P \; \text{rough}}} \frac{b(N)}{N^{\sigma+it}} \Biggr|^2 dt , \nonumber
\end{eqnarray}
where we temporarily set $b(N) := \#\{(n_1,n_2): n_i \leq T^{1/2 - 2\epsilon} \; \text{are} \; P \; \text{rough}, \; n_1 n_2 = N\}$. By Lemma \ref{meanvaluelem}, this is all
$$ \ll (\frac{\log P}{V})^2 \cdot T \sum_{\substack{m \leq T^{\epsilon}, \\ m \; \text{is} \; P \; \text{smooth}}} \frac{1}{m} \sum_{\substack{N \leq T^{1 - 4\epsilon}, \\ N \; \text{is} \; P \; \text{rough}}} \frac{b(N)^2}{N^{2\sigma}} \ll (\frac{\log P}{V})^2 \cdot T \log P \sum_{\substack{N \leq T^{1 - 4\epsilon}, \\ N \; \text{is} \; P \; \text{rough}}} \frac{d(N)^2}{N^{2\sigma}} , $$
where $d(N)$ is the classical divisor function and we again used the standard estimate $\sum_{\substack{m \leq T^{\epsilon}, \\ m \; \text{is} \; P \; \text{smooth}}} \frac{1}{m} \leq \prod_{p \leq P} (1 - \frac{1}{p})^{-1} \ll \log P$. Finally, if we upper bound the sum over $N$ by the infinite series $\sum_{\substack{N : \\ p | N \Rightarrow P < p \leq T}} \frac{d(N)^2}{N^{2\sigma}} \leq \prod_{P < p \leq T} (1 - \frac{4}{p^{2\sigma}})^{-1} \ll (\frac{\log T}{\log P})^4$, the lemma follows.
\end{proof}

We shall also require an extra probabilistic estimate reflecting the fact that we are looking for very large values, to feed into the general machinery from the author's paper~\cite{harperrmflow} in place of the upper barrier estimate that was originally used there.

\begin{probres1}[Version of the Ballot Theorem]
Let $a, b$ and $n \in \N$ be large (i.e. larger than certain absolute positive constants), and let $G_1 , ..., G_n$ be independent Gaussian random variables, each having mean zero and variance between $1/20$ and $20$ (say). Then we have the uniform upper bound
$$ \p(\sum_{m=1}^{j} G_m \leq a \; \forall 1 \leq j \leq n, \; \text{and} \; a - b \leq \sum_{m=1}^{n} G_m \leq a ) \ll \min\{1, \frac{a}{\sqrt{n}}\} (\min\{1, \frac{b}{\sqrt{n}}\})^2 . $$
\end{probres1}

Results of approximately this form, under various different assumptions on the $G_i$ and on $a,b,n$, are fairly standard. Since it is very neat and rather short, we shall give a full proof for our situation roughly following an argument from section 6 of Webb~\cite{webb}.

\begin{proof}[Proof of Probability Result 1]
Note that if the event in the result occurs, then in particular at the end of our random walk we must have
$$ \sum_{m=j}^{n} G_m = \sum_{m=1}^{n} G_m - \sum_{m=1}^{j-1} G_m \geq (a-b) - a = - b \;\;\; \forall 2n/3 < j \leq n , $$
say. We also must have $\sum_{m=1}^{j} G_m \leq a$ for all $j \leq n/3$, and we must have
$$ \sum_{n/3 < m \leq 2n/3} G_m \in [a - b - \sum_{m \notin (n/3, 2n/3]} G_m, a - \sum_{m \notin (n/3, 2n/3]} G_m] . $$

Now $\sum_{n/3 < m \leq 2n/3} G_m$ is a Gaussian random variable, with mean zero and variance of order $n$, that is independent of $(G_m)_{m \notin (n/3, 2n/3]}$. The probability that such a Gaussian takes values in an interval of length $b$ is $\ll \min\{1, \frac{b}{\sqrt{n}}\}$, so the probability of the event in the result is
$$ \ll \min\{1, \frac{b}{\sqrt{n}}\} \cdot \p(\sum_{m=1}^{j} G_m \leq a \; \forall 1 \leq j \leq n/3, \; \text{and} \; - ( \sum_{m=j}^{n} G_m) \leq b \; \forall 2n/3 < j \leq n) . $$
By independence of the $G_m$, this is all equal to
$$ \min\{1, \frac{b}{\sqrt{n}}\} \cdot \p(\sum_{m=1}^{j} G_m \leq a \; \forall 1 \leq j \leq n/3) \cdot \p( - ( \sum_{m=j}^{n} G_m) \leq b \; \forall 2n/3 < j \leq n) , $$
and using standard results on the maximum of random walks (see e.g. estimate (A1) from the appendix of Harper~\cite{harperrmflow}) we have $\p(\sum_{m=1}^{j} G_m \leq a \; \forall 1 \leq j \leq n/3) \asymp \min\{1, \frac{a}{\sqrt{n}}\}$ and $\p( - ( \sum_{m=j}^{n} G_m) \leq b \; \forall 2n/3 < j \leq n) \asymp \min\{1, \frac{b}{\sqrt{n}}\}$, which yields the claimed result.
\end{proof}

It isn't too hard to see that provided $a \ll \sqrt{n}$, Probability Result 1 is sharp (because we may assume for a lower bound that $b \leq \sqrt{n}$ as well, and roughly speaking then the three bounds used in the above proof are sharp and, with positive conditional probability, if all three events hold then the random walk will behave in the way required by the Ballot Theorem). But we won't need to know this for the proof of the upper bound in Theorem 2.

\section{Proof of Theorem 2}
Many of the arguments we used previously can be transferred over to handle Theorem 2, but for one technical reason later in the proof (see Lemma \ref{eulerprodapproxlem} below) we shall slightly change our definition of the approximating points $h(j)$ (forcing them to be closer to $h$), and related objects. So for $P$ a large parameter and $|h| \leq 1/2$, let us define $\tilde{h}(-1) := h$, and then for $0 \leq j \leq \log\log P - 1$ set
$$ \tilde{h}(j) := \max\{u \leq \tilde{h}(j-1): u = \frac{n}{((\log P)/e^j) (\log\log P)^3} \; \text{for some} \; n \in \Z\} . $$
Also let $\tilde{h}(*)$ denote the point of the form $n/\log T$ that is closest to $h$. We let $ I_l(h) = I_{l,t}(h)$ denote the partial Euler products, exactly as in section \ref{thm1proofsec}.

Next, for each $|h| \leq 1/2$ let $\tilde{\mathcal{G}}(h) = \tilde{\mathcal{G}}_{t}(h)$ denote the ``good'' event that for all $0 \leq j \leq \log\log P - B - 1$, we have
\begin{equation}\label{tilgdef}
\Biggl|\sum_{l=j}^{\lfloor \log\log P \rfloor - B - 1} \log I_{l,t}(\tilde{h}(l)) \Biggr| \leq \log\log P - j + 3\log\log\log P + U .
\end{equation}
Here $B \in \N$ will be a certain large constant as before, and $U \geq 0$ is a parameter. We also let $\tilde{\mathcal{G}} = \tilde{\mathcal{G}}_{t}$ denote the event that $\tilde{\mathcal{G}}_{t}(h)$ holds simultaneously for all $|h| \leq 1/2$. Note this is all very similar to our set-up with $\mathcal{G}_{t}$ from the proof of Theorem 1, except that \eqref{tilgdef} involves $I_{l,t}(\tilde{h}(l))$ rather than its absolute value $|I_{l,t}(\tilde{h}(l))|$. This small change will be useful later (see the proof of Key Proposition 4, below) when we come to compare certain sums with (essentially) Euler products, and make no difference to our analysis of how frequently $\tilde{\mathcal{G}}_{t}$ holds.

\vspace{12pt}
Our key estimates now shall be the following.

\begin{keyprop3}
Uniformly for all large $P \leq T^{1/(\log\log T)^2}$ and $0 \leq U \leq 2\log\log P$, say, we have
$$ \frac{1}{T} \text{meas}\{T \leq t \leq 2T : \tilde{\mathcal{G}}_{t} \; \text{fails}\} \ll e^{-2U} . $$
\end{keyprop3}

It is perhaps worth noting that the range of $U$ allowed in Key Proposition 3 could easily be increased at the cost of a stronger restriction on $P$, and the bound can also be strengthened when $P$ is large (one gets an additional factor of the shape $e^{-\Theta(U^{2}/\log\log P)}$ multiplying the right hand side). But the stated version will suffice for our purposes. 

\begin{keyprop4}
Uniformly for all large $P \leq T^{1/(\log\log T)^{8}}$ and $\frac{20\log P \log\log P}{\log T} < \epsilon < 1/10$ and $0 \leq U \leq 2\log\log P$ and $V \geq e^{-U}$ (say); and for all $|h| \leq 1/2$ and $\frac{1}{2} - \frac{1}{\log T} \leq \sigma \leq \frac{1}{2} + \frac{1}{\log T}$; we have
\begin{eqnarray}
&& \frac{1}{T} \int_{T}^{2T} \textbf{1}_{\tilde{\mathcal{G}}_{t}(h)} \textbf{1}_{|\sum_{\substack{m \leq T^{\epsilon}, \\ P \; \text{smooth}}} \frac{1}{m^{1/2+i(t+h)}}| > (\log P)/V} |\sum_{\substack{m \leq T^{\epsilon}, \\ P \; \text{smooth}}} \frac{1}{m^{1/2+i(t+h)}}|^2 |\sum_{\substack{n \leq T^{1/2 - 2\epsilon}, \\ n \; \text{is} \; P \; \text{rough}}} \frac{1}{n^{\sigma+i(t+h)}}|^2 dt \nonumber \\
& \ll & \log T \min\{1, \frac{\log\log\log P + U}{\sqrt{\log\log P}} \} \min\{1, \frac{\log\log\log P + U + \log V}{\sqrt{\log\log P}} \}^2 , \nonumber
\end{eqnarray}
where $\textbf{1}$ denotes the indicator function.
\end{keyprop4}

\begin{proof}[Proof of Theorem 2, assuming Key Propositions 3 and 4]
We again fix the choices $\epsilon = \frac{1}{(\log\log T)^2}$ and $P = T^{1/(\log\log T)^{8}}$, as well as setting $V = e^{-U} (\log\log T)^6$. Recall that in Theorem 2 we have $0 \leq U \leq \log\log T$, so all these values satisfy the conditions of the Key Propositions.

For each $T \leq t \leq 2T$, let $h_t$ denote a value of $|h| \leq 1/2$ at which $\max_{|h| \leq 1/2} |\zeta(1/2+it+ih)|$ is attained. (This notation will be temporary, so shouldn't cause confusion with our existing notation for the approximating points $\tilde{h}(l)$.) Thus our goal is to show that
$$ \frac{1}{T} \text{meas}\{T \leq t \leq 2T : |\zeta(1/2+it+ih_t)| \geq \frac{e^{U} \log T}{(\log\log T)^{3/4}} \} \ll e^{-2U} (\log\log\log T + U) (\log\log\log T)^2 . $$
By Key Proposition 3, the measure of the set of $t$ for which $\tilde{\mathcal{G}}_{t}$ fails is smaller than this bound, so where helpful we may restrict attention to those $t$ for which $\tilde{\mathcal{G}}_{t}$ holds.

By the approximate functional equation (Zeta Function Result 1), if $|\zeta(1/2+it+ih_t)| \geq \frac{e^{U} \log T}{(\log\log T)^{3/4}}$ then we must have $\sum_{n \leq \sqrt{(t + h_t)/2\pi}} \frac{1}{n^{1/2 + it + ih_t}} \gg \frac{e^{U} \log T}{(\log\log T)^{3/4}}$. Then applying Cauchy's Integral Formula to the holomorphic function $s \mapsto \sum_{n \leq \sqrt{(t + h_t)/2\pi}} \frac{1}{n^{s}}$, we have
$$ \sum_{n \leq \sqrt{(t + h_t)/2\pi}} \frac{1}{n^{1/2 + it + ih_t}} = \frac{1}{2\pi i} \left( \int_{A_t} + \int_{B_t} + \int_{C_t} + \int_{D_t} \right) \frac{ \sum_{n \leq \sqrt{(t + h_t)/2\pi}} \frac{1}{n^s}}{s - (1/2 + it + ih_t)} ds , $$
where $A_t , B_t , C_t , D_t$ are the sides (taken with anticlockwise orientation) of the small box whose vertices have real part $1/2 \pm 1/\log T$ and imaginary part $t + \tilde{h_t}(*) \pm \frac{1}{\log T}$. Note that by definition of $\tilde{h_t}(*)$ we have $|h_t - \tilde{h_t}(*)| \leq 1/(2\log T)$, and therefore not only is the point $1/2 + it + ih_t$ certainly inside the small box, but we also have $|s - (1/2 + it + ih_t)| \asymp 1/\log T$ on all four sides. Note also that we may replace $\sum_{n \leq \sqrt{(t + h_t)/2\pi}} \frac{1}{n^s}$ in all four integrals by $\sum_{n \leq \sqrt{\Im(s)/2\pi}} \frac{1}{n^s}$, at the cost of a tiny error term of size at most $O(T^{-1/4})$. So if $\sum_{n \leq \sqrt{(t + h_t)/2\pi}} \frac{1}{n^{1/2 + it + ih_t}} \gg \frac{e^{U} \log T}{(\log\log T)^{3/4}}$, then we must have $\int |\sum_{n \leq \sqrt{\Im(s)/2\pi}} \frac{1}{n^s}| d|s| \gg \frac{e^{U}}{(\log\log T)^{3/4}}$ for at least one of the four integrals $\int_{A_t}, \int_{B_t}, \int_{C_t}, \int_{D_t}$.

First we shall consider what happens for the integrals over the two vertical sides. The treatment of both is exactly the same, so to simplify the writing let $\sigma$ denote either $1/2 + 1/\log T$ or $1/2 - 1/\log T$. We have
\begin{eqnarray}\label{cidecomp}
&& \int_{\tilde{h_t}(*) - \frac{1}{\log T}}^{\tilde{h_t}(*) + \frac{1}{\log T}} \Biggl|\sum_{n \leq \sqrt{(t+h)/2\pi}} \frac{1}{n^{\sigma + it+ih}} \Biggr| dh \\
& \leq & \int_{\tilde{h_t}(*) - \frac{1}{\log T}}^{\tilde{h_t}(*) + \frac{1}{\log T}} \Biggl|\sum_{n \leq \sqrt{(t+h)/2\pi}} \frac{1}{n^{\sigma + it+ih}} - \sum_{\substack{m \leq T^{\epsilon}, \\ m \; \text{is} \; P \; \text{smooth}}} \frac{1}{m^{1/2+it+ih}} \sum_{\substack{n \leq T^{1/2 - 2\epsilon}, \\ n \; \text{is} \; P \; \text{rough}}} \frac{1}{n^{\sigma+it+ih}} \Biggr| dh \nonumber \\
&& + \int_{\tilde{h_t}(*) - \frac{1}{\log T}}^{\tilde{h_t}(*) + \frac{1}{\log T}} \textbf{1}_{|\sum_{\substack{m \leq T^{\epsilon}, \\ P \; \text{smooth}}} \frac{1}{m^{1/2+it+ih}}| \leq (\log P)/V} \Biggl|\sum_{\substack{m \leq T^{\epsilon}, \\ P \; \text{smooth}}} \frac{1}{m^{1/2+it+ih}} \sum_{\substack{n \leq T^{1/2 - 2\epsilon}, \\ n \; \text{is} \; P \; \text{rough}}} \frac{1}{n^{\sigma+it+ih}} \Biggr| dh \nonumber \\
&& + \int_{\tilde{h_t}(*) - \frac{1}{\log T}}^{\tilde{h_t}(*) + \frac{1}{\log T}} \textbf{1}_{|\sum_{\substack{m \leq T^{\epsilon}, \\ P \; \text{smooth}}} \frac{1}{m^{1/2+it+ih}}| > (\log P)/V} \Biggl|\sum_{\substack{m \leq T^{\epsilon}, \\ P \; \text{smooth}}} \frac{1}{m^{1/2+it+ih}} \sum_{\substack{n \leq T^{1/2 - 2\epsilon}, \\ n \; \text{is} \; P \; \text{rough}}} \frac{1}{n^{\sigma+it+ih}} \Biggr| dh , \nonumber
\end{eqnarray}
so if the left hand side is $\gg \frac{e^{U}}{(\log\log T)^{3/4}}$ then the same must be true for at least one of the three terms on the right. By the Cauchy--Schwarz inequality, if the first term is $\gg \frac{e^{U}}{(\log\log T)^{3/4}}$ then we must have
$$ (\frac{e^{U}}{(\log\log T)^{3/4}})^2 \ll \frac{1}{\log T} \int_{\tilde{h_t}(*) - \frac{1}{\log T}}^{\tilde{h_t}(*) + \frac{1}{\log T}} \Biggl|\sum_{n \leq \sqrt{\frac{t+h}{2\pi}}} \frac{1}{n^{\sigma + i(t+h)}} - \sum_{\substack{m \leq T^{\epsilon}, \\ P \; \text{smooth}}} \frac{1}{m^{1/2+i(t+h)}} \sum_{\substack{n \leq T^{1/2 - 2\epsilon}, \\ n \; \text{is} \; P \; \text{rough}}} \frac{1}{n^{\sigma+i(t+h)}} \Biggr|^2 . $$
Having efficiently applied the Cauchy--Schwarz inequality to this short integral, we can now upper bound the right hand side by increasing the range of integration to all $|h| \leq 1/2$ (thereby removing the dependence on $t$ there). Then integrating over $T \leq t \leq 2T$ and applying Lemma \ref{polyeditlem}, we deduce that the measure of the set of $t$ for which this inequality holds is
$$ \ll (\frac{(\log\log T)^{3/4}}{e^U})^2 T (e^{-\epsilon (\log T)/\log P} + \epsilon) \ll \frac{e^{-2U} T}{\sqrt{\log\log T}} , $$
which is more than good enough. Similarly, if the second of our three integrals is $\gg \frac{e^{U}}{(\log\log T)^{3/4}}$ then, applying H\"{o}lder's inequality before again increasing the range of integration to all $|h| \leq 1/2$, we must have
$$ (\frac{e^{U}}{(\log\log T)^{3/4}})^4 \ll \frac{1}{\log^{3}T} \int_{-\frac{1}{2}}^{\frac{1}{2}} \textbf{1}_{|\sum_{\substack{m \leq T^{\epsilon}, \\ P \; \text{smooth}}} \frac{1}{m^{1/2+i(t+h)}}| \leq \frac{\log P}{V}} \Biggl|\sum_{\substack{m \leq T^{\epsilon}, \\ P \; \text{smooth}}} \frac{1}{m^{1/2+i(t+h)}} \sum_{\substack{n \leq T^{1/2 - 2\epsilon}, \\ P \; \text{rough}}} \frac{1}{n^{\sigma+i(t+h)}} \Biggr|^4 dh . $$
Integrating over $T \leq t \leq 2T$ and applying Lemma \ref{fourthmomentlem}, we find the measure of the set of $t$ for which this holds is
$$ \ll (\frac{(\log\log T)^{3/4}}{e^U})^4 \frac{T}{V^2} \frac{\log T}{\log P} \ll \frac{e^{-4U} T (\log\log T)^{11}}{V^2} \ll \frac{e^{-2U} T}{\log\log T} , $$
which again is more than good enough.

Finally we turn to the third integral from \eqref{cidecomp}. If this is $\gg \frac{e^{U}}{(\log\log T)^{3/4}}$ then, applying the Cauchy--Schwarz inequality before increasing the range of integration, we must have
$$ (\frac{e^{U}}{(\log\log T)^{3/4}})^2 \ll \frac{1}{\log T} \int_{-\frac{1}{2}}^{\frac{1}{2}} \textbf{1}_{|\sum_{\substack{m \leq T^{\epsilon}, \\ P \; \text{smooth}}} \frac{1}{m^{1/2+i(t+h)}}| > \frac{\log P}{V}} \Biggl|\sum_{\substack{m \leq T^{\epsilon}, \\ P \; \text{smooth}}} \frac{1}{m^{1/2+i(t+h)}} \sum_{\substack{n \leq T^{1/2 - 2\epsilon}, \\ P \; \text{rough}}} \frac{1}{n^{\sigma+i(t+h)}} \Biggr|^2 dh . $$
Then the measure of the set of $T \leq t \leq 2T$ for which this inequality holds, {\em and for which $\tilde{\mathcal{G}}_{t}$ also holds}, is bounded by a constant times
$$ \frac{(\log\log T)^{3/2}}{e^{2U} \log T} \int_{-\frac{1}{2}}^{\frac{1}{2}} \int_{T}^{2T} \textbf{1}_{\tilde{\mathcal{G}}_{t}} \textbf{1}_{|\sum_{\substack{m \leq T^{\epsilon}, \\ P \; \text{smooth}}} \frac{1}{m^{1/2+it+ih}}| > \frac{\log P}{V}} \Biggl|\sum_{\substack{m \leq T^{\epsilon}, \\ P \; \text{smooth}}} \frac{1}{m^{1/2+it+ih}} \sum_{\substack{n \leq T^{1/2 - 2\epsilon}, \\ P \; \text{rough}}} \frac{1}{n^{\sigma+it+ih}} \Biggr|^2 dt dh . $$
In the inner integral we can now upper bound $\textbf{1}_{\tilde{\mathcal{G}}_{t}}$ by $\textbf{1}_{\tilde{\mathcal{G}}_{t}(h)}$, and then Key Proposition 4 implies this is all
\begin{eqnarray}
& \ll & \frac{(\log\log T)^{3/2}}{e^{2U}} T \min\{1, \frac{\log\log\log P + U}{\sqrt{\log\log P}} \} \min\{1, \frac{\log\log\log P + U + \log V}{\sqrt{\log\log P}} \}^2 \nonumber \\
& \ll & e^{-2U} T (\log\log\log T + U) (\log\log\log T)^2 , \nonumber
\end{eqnarray}
as required.

It only remains to deal with the integrals over the horizontal sides of our original box, which are of the form
$$ \int_{1/2 - 1/\log T}^{1/2 + 1/\log T} \Biggl|\sum_{n \leq \sqrt{\frac{t+\tilde{h_t}(*) - 1/\log T}{2\pi}}} \frac{1}{n^{\sigma + i(t+ \tilde{h_t}(*) - 1/\log T)}} \Biggr| d\sigma , $$
and the same with $t+ \tilde{h_t}(*) - 1/\log T$ replaced by $t+ \tilde{h_t}(*) + 1/\log T$. We can upper bound this by three integrals in the same way as for \eqref{cidecomp}, and now thanks to the definition of $\tilde{h_t}(*)$, the value $\tilde{h_t}(*) \pm 1/\log T$ is always of the form $n/\log T$ for some {\em integer} $|n| \leq (\log T)/2 + O(1)$ (i.e. it always belongs to a fixed set of $\asymp \log T$ values). So after applying the Cauchy--Schwarz inequality or H\"{o}lder's inequality to the integrals as before, we can remove the dependence on $t$ coming from $\tilde{h_t}(*)$ by summing over all possible values $n/\log T$ of the shift $\tilde{h_t}(*) \pm 1/\log T$ (rather than by changing the range of integration over $\sigma$). The remainder of the argument goes through without change from the case of the vertical sides.
\end{proof}

\vspace{12pt}
Once again, it remains to prove the Key Propositions.

\begin{proof}[Proof of Key Proposition 3]
The proof proceeds using the union bound, exactly as in the proof of Key Proposition 2. Note that we don't need to make any special changes to account for the fact that we are interested in $\log I_{l,t}(\tilde{h}(l))$ rather than $\log|I_{l,t}(\tilde{h}(l))| = \Re \log I_{l,t}(\tilde{h}(l))$, since in the proof of Key Proposition 2 we already upper bound the real parts of the prime number sums by their absolute values.

We give a few more details. We need to bound the measure of the set of $T \leq t \leq 2T$ for which
$$ \Biggl|\sum_{l=j}^{\lfloor \log\log P \rfloor - B - 1} \sum_{\substack{P^{e^{-(l+1)}} < p \leq P^{e^{-l}}, \\ v=1,2}} \frac{1}{vp^{v(1/2+it+i\tilde{h}(l))}} \Biggr| > \log\log P - j + 3\log\log\log P + U + O(1) . $$
As in the proof of Key Proposition 2, for any $k \in \N$ such that $P^{2k e^{-j}} < T$ we can upper bound this measure by
$$ \ll T \sqrt{k} \left(\frac{k(\log\log P - j + O(1))}{e(\log\log P - j + 3\log\log\log P + U + O(1))^{2}} \right)^{k} . $$
Again, choosing $k = \lfloor (\log\log P - j + 3\log\log\log P + U + O(1))^{2}/(\log\log P - j + O(1)) \rfloor$ is roughly optimal, and yields an overall bound
$$ \ll T \sqrt{(\log\log P - j)} \frac{e^j}{\log P} e^{-2(3 \log\log\log P + U)} \ll T \frac{e^j}{\log P} e^{-2U} \frac{1}{(\log\log P)^{11/2}} . $$
Since we assume that $U \leq 2\log\log P$ we have $k \leq (1/3)(\log\log P)^2$ here (provided the constant $B$ is fixed sufficiently large so the denominator $\log\log P - j + O(1) \geq B + O(1)$ is always large). Thus the condition $P^{2k e^{-j}} < T$ will certainly be satisfied, because of our assumption that $P \leq T^{1/(\log\log T)^2}$.

The proof ends in the same way as for Key Proposition 2, by summing up over all $\tilde{h}(j)$ (of which there are now $\ll \frac{\log P}{e^j} (\log\log P)^3$ values) and over all $0 \leq j \leq \log\log P - B - 1$.
\end{proof}

The proof of Key Proposition 4 will be more difficult, although several of the technical issues have already been dealt with in the proof of Key Proposition 1. We shall also need a couple of lemmas. The first shows that, in an average sense, we can compare and replace $\sum_{\substack{m \leq T^{\epsilon}, \\ m \; \text{is} \; P \; \text{smooth}}} \frac{1}{m^{1/2+i(t+h)}}$ with something more like an Euler product. Using the fact that the approximating points $\tilde{h}(l)$ are all now a bit closer to $h$, we can actually replace $\sum_{\substack{m \leq T^{\epsilon}, \\ m \; \text{is} \; P \; \text{smooth}}} \frac{1}{m^{1/2+i(t+h)}}$ by an Euler product-like object with $h$ replaced by the $\tilde{h}(l)$.

\begin{lem6}\label{eulerprodapproxlem}
Let $P \leq T^{1/(100 \log\log T)}$ be large, and let $\frac{e^2 \log P \log\log P}{\log T} < \epsilon < 1/10$. For each prime $p \leq P$, define $l(p) \in \N \cup \{0\}$ to be the value of $l$ for which $P^{e^{-(l+1)}} < p \leq P^{e^{-l}}$. Then uniformly for all $|h| \leq 1/2$ and $\frac{1}{2} - \frac{1}{\log T} \leq \sigma \leq \frac{1}{2} + \frac{1}{\log T}$, we have
\begin{eqnarray}
&& \int_{T}^{2T} \Biggl|\sum_{\substack{m \leq T^{\epsilon}, \\ P \; \text{smooth}}} \frac{1}{m^{1/2+i(t+h)}} - \sum_{k=0}^{\lfloor \frac{\epsilon \log T}{\log P} \rfloor} \frac{1}{k!} (\sum_{p^{j} \leq P} \frac{1}{j p^{j(1/2 + it + i\tilde{h}(l(p)))}} )^k \Biggr|^2 \Biggl|\sum_{\substack{n \leq T^{1/2 - 2\epsilon}, \\ P \; \text{rough}}} \frac{1}{n^{\sigma+i(t+h)}} \Biggr|^2 dt \nonumber \\
& \ll & \frac{T \log T}{(\log\log P)^2} . \nonumber
\end{eqnarray}
\end{lem6}

\begin{proof}[Proof of Lemma \ref{eulerprodapproxlem}]
If we expand out $\sum_{k=0}^{\lfloor \frac{\epsilon \log T}{\log P} \rfloor} \frac{1}{k!} (\sum_{p^{j} \leq P} \frac{1}{j p^{j(1/2 + it)}} )^k$, we obtain a Dirichlet polynomial of the form $\sum_{\substack{m \leq T^{\epsilon}, \\ m \; \text{is} \; P \; \text{smooth}}} \frac{c(m)}{m^{1/2+it}}$, for certain coefficients $c(m)$ that we shall investigate below. Noting this, and using the second part of Lemma \ref{meanvaluelem}, we deduce that the left hand side in Lemma \ref{eulerprodapproxlem} is
\begin{eqnarray}
& \ll & T \sum_{\substack{n \leq T^{1/2 - 2\epsilon}, \\ n \; \text{is} \; P \; \text{rough}}} \frac{1}{n^{2\sigma}} \sum_{\substack{m \leq T^{\epsilon}, \\ m \; \text{is} \; P \; \text{smooth}}} \frac{|m^{-ih} - c(m) \prod_{p^j || m} p^{-ij\tilde{h}(l(p))}|^2}{m} \nonumber \\
& \ll & T \frac{\log T}{\log P} \sum_{\substack{m \leq T^{\epsilon}, \\ m \; \text{is} \; P \; \text{smooth}}} \frac{|1 - c(m) \prod_{p^j || m} p^{-ij(\tilde{h}(l(p)) - h)}|^2}{m} , \nonumber
\end{eqnarray}
where $\prod_{p^j || m}$ denotes the product over all the largest powers of primes that divide $m$ (so that $m = \prod_{p^j || m} p^j$).

Next, using the Taylor expansions of the exponential and the logarithm we have
$$  \sum_{\substack{m : \\ P \; \text{smooth}}} \frac{1}{m^{1/2+it}} = \prod_{p \leq P} (1 - \frac{1}{p^{1/2+it}})^{-1} = \exp\{-\sum_{p \leq P} \log(1 - \frac{1}{p^{1/2+it}})\} = \sum_{k=0}^{\infty} \frac{1}{k!} (\sum_{p^{j} : p \leq P} \frac{1}{j p^{j(1/2 + it)}} )^k , $$
so comparing with the truncated series $\sum_{k=0}^{\lfloor \frac{\epsilon \log T}{\log P} \rfloor} \frac{1}{k!} (\sum_{p^{j} \leq P} \frac{1}{j p^{j(1/2 + it)}} )^k$ we see that $0 \leq c(m) \leq 1$. Thus we can upper bound our sum over $m \leq T^{\epsilon}$ by
\begin{eqnarray}
& \ll & \sum_{\substack{m \leq T^{\epsilon}, \\ m \; \text{is} \; P \; \text{smooth}}} \frac{|1 - c(m)|^2}{m} + \sum_{\substack{m \leq T^{\epsilon}, \\ m \; \text{is} \; P \; \text{smooth}}} \frac{|1 - \prod_{p^j || m} p^{-ij(\tilde{h}(l(p)) - h)}|^2}{m} \nonumber \\
& \ll &  \sum_{\substack{m \leq T^{\epsilon}, \\ m \; \text{is} \; P \; \text{smooth}}} \frac{1 - c(m)}{m} + \sum_{\substack{m \leq T^{\epsilon}, \\ m \; \text{is} \; P \; \text{smooth}}} \frac{\sum_{p^j || m} j (\log p) |\tilde{h}(l(p)) - h|}{m} . \nonumber
\end{eqnarray}
Swapping the order of the summations, the second term here is
$$ \leq \sum_{\substack{p^j \leq T^{\epsilon} : \\ p \leq P}} \frac{j (\log p) |\tilde{h}(l(p)) - h|}{p^j} \sum_{\substack{n \leq T^{\epsilon}/p^j , \\ n \; \text{is} \; P \; \text{smooth}}} \frac{1}{n} \ll \log P \sum_{p \leq P} \frac{(\log p) |\tilde{h}(l(p)) - h|}{p} . $$
And by definition of $l(p)$ and of the approximating points $\tilde{h}$, we always have $|\tilde{h}(l(p)) - h| \ll \frac{e^{l(p)}}{\log P (\log\log P)^3} \ll \frac{1}{\log p (\log\log P)^3}$, and so this sum over $p$ is $\ll \sum_{p \leq P} \frac{1}{p (\log\log P)^3} \ll \frac{1}{(\log\log P)^2}$ by the Mertens estimate. This contribution is acceptable for the lemma.

Finally, by definition of $c(m)$ the remaining first sum $\sum_{\substack{m \leq T^{\epsilon}, \\ P \; \text{smooth}}} \frac{1 - c(m)}{m}$ is at most
$$ \sum_{\substack{m : \\ P \; \text{smooth}}} \frac{1}{m} - \sum_{k=0}^{\lfloor \frac{\epsilon \log T}{\log P} \rfloor} \frac{1}{k!} (\sum_{p^{j} \leq P} \frac{1}{j p^{j}} )^k = \exp\{ \sum_{p^{j} \leq P} \frac{1}{j p^{j}} + O(\frac{1}{\sqrt{P} \log P})\} - \sum_{k=0}^{\lfloor \frac{\epsilon \log T}{\log P} \rfloor} \frac{1}{k!} (\sum_{p^{j} \leq P} \frac{1}{j p^{j}} )^k , $$
since $\sum_{\substack{m : \\ P \; \text{smooth}}} \frac{1}{m} = \prod_{p \leq P} (1 - \frac{1}{p})^{-1} = \exp\{ \sum_{p^{j} : p \leq P} \frac{1}{j p^{j}}\}$, and those terms in the exponential with $p \leq P$ but $p^j > P$ contribute at most $O(\frac{1}{\sqrt{P} \log P})$. Since $\sum_{p^{j} \leq P} \frac{1}{j p^{j}} = \log\log P + O(1)$ by the Mertens estimate, this is all equal to
$$ \exp\{ \sum_{p^{j} \leq P} \frac{1}{j p^{j}}\} + O(\frac{1}{\sqrt{P}}) - \sum_{k=0}^{\lfloor \frac{\epsilon \log T}{\log P} \rfloor} \frac{1}{k!} (\sum_{p^{j} \leq P} \frac{1}{j p^{j}} )^k = \sum_{k= \lfloor \frac{\epsilon \log T}{\log P} \rfloor + 1}^{\infty} \frac{1}{k!} (\sum_{p^{j} \leq P} \frac{1}{j p^{j}} )^k + O(\frac{1}{\sqrt{P}}) . $$
And using our assumption that $\frac{\epsilon \log T}{\log P} \geq e^2 \log\log P$, we see the sum over $k$ is dominated by the term with $k= \lfloor \frac{\epsilon \log T}{\log P} \rfloor + 1$. Inserting the lower bound $k! \geq (k/e)^k$, a quick calculation shows that term is $\ll e^{-\epsilon (\log T)/\log P}$, which is more than good enough for the lemma.
\end{proof}

Similarly as for Key Proposition 1, in the proof of Key Proposition 4 we will ultimately reduce matters to calculating an average involving random multiplicative functions. The following lemma, which combines Probability Result 1 with some machinery from the author's paper~\cite{harperrmflow}, will provide a suitable estimate for that average. For each $|h| \leq 1/2$, let $\tilde{\mathcal{G}}'(h)$ denote the event that: for all $0 \leq j \leq \log\log P - B - 1$, we have
$$ \Biggl( \frac{\log P}{e^j} e^{3\log\log\log P + U + 19/3} \Biggr)^{-1} \leq \prod_{l = j}^{\lfloor \log\log P \rfloor - B - 1} |I_{l, \text{rand}}(\tilde{h}(l))| \leq \frac{\log P}{e^j} e^{3\log\log\log P + U + 19/3} , $$
and also that for $j=0$ we have the (generally more stringent) lower bound
$$ \prod_{l = 0}^{\lfloor \log\log P \rfloor - B - 1} |I_{l, \text{rand}}(\tilde{h}(l))| \geq \frac{\log P}{V B_0 e^{19/3}} . $$
Here $I_{l, \text{rand}}(\cdot)$ and $\tilde{h}(l)$ are the random Euler products and approximating points that we introduced earlier, and $B_0$ denotes a certain large absolute constant coming from the proof of Key Proposition 4 (see below).

\begin{lem7}\label{mainprob2lem}
Let $f(n)$ denote a Steinhaus random multiplicative function. Uniformly for all large $P \leq \sqrt{T}$ and $0 < \epsilon < 1/10$ and $0 \leq U \leq 2\log\log P$ (say) and $V \geq e^{-U}$; and for all $|h| \leq 1/2$ and $\frac{1}{2} - \frac{1}{\log T} \leq \sigma \leq \frac{1}{2} + \frac{1}{\log T}$; we have
\begin{eqnarray}
&& \E \textbf{1}_{\tilde{\mathcal{G}}^{'}(h)} |\sum_{\substack{m \leq T^{\epsilon}, \\ m \; \text{is} \; P \; \text{smooth}}} \frac{f(m)}{m^{1/2+ih}}|^2 |\sum_{\substack{n \leq T^{1/2 - 2\epsilon}, \\ n \; \text{is} \; P \; \text{rough}}} \frac{f(n)}{n^{\sigma+ih}}|^2 \nonumber \\
& \ll & \log T \left( \min\{1, \frac{\log\log\log P + U}{\sqrt{\log\log P}} \} \min\{1, \frac{\log\log\log P + U + \log V}{\sqrt{\log\log P}} \}^2 + e^{-\epsilon (\log T)/\log P} \right) . \nonumber
\end{eqnarray}
\end{lem7}

Once again, the restriction that $U \leq 2\log\log P$ is not very essential here, but will slightly simplify the writing of the proof by allowing us to perform some truncation steps without explicit reference to $U$, only in terms of $P$. A condition roughly of the form $V \geq e^{-U}$ is rather natural, to ensure that $\log\log\log P + U + \log V$ is large (and if $V$ were significantly smaller, the event $\tilde{\mathcal{G}}^{'}(h)$ could in fact never occur because the lower and upper bounds for $j=0$ would be incompatible).

\begin{proof}[Proof of Lemma \ref{mainprob2lem}]
Arguing exactly as in the proof of Lemma \ref{mainprob1lem} above, we can reduce the problem to showing that
$$ \E \textbf{1}_{\tilde{\mathcal{G}}^{'}(h)} |\prod_{p \leq P} (1 - \frac{f(p)}{p^{1/2+ih}})|^{-2} \ll \log P \min\{1, \frac{\log\log\log P + U}{\sqrt{\log\log P}} \} \min\{1, \frac{\log\log\log P + U + \log V}{\sqrt{\log\log P}} \}^2 . $$
And we may assume that $V \leq e^{\sqrt{\log\log P}}$, otherwise we can just discard the extra lower bound condition on $\prod_{l = 0}^{\lfloor \log\log P \rfloor - B - 1} |I_{l, \text{rand}}(\tilde{h}(l))|$ from the definition of $\tilde{\mathcal{G}}^{'}(h)$, and conclude as in the proof of Lemma \ref{mainprob1lem}.

Now our desired upper bound almost follows just by combining Probability Result 1 with Lemma 4 of Harper~\cite{harperrmflow}, except that very large values of $l$ in the products (corresponding to fairly small primes) create a technical problem with the conditions of that lemma. We will describe how to circumvent this problem. The quantity $13\log\log\log P + U$ will occur a few times, so we shall temporarily write $E(P,U) := 13\log\log\log P + U$. Indeed, note that if $\tilde{\mathcal{G}}^{'}(h)$ holds, and if we temporarily set $\tilde{F} = \tilde{F}(h) := \prod_{l = \lfloor \log\log P - 10\log\log\log P \rfloor}^{\lfloor \log\log P \rfloor - B - 1} |I_{l, \text{rand}}(\tilde{h}(l))|^{-1}$, then we must have
$$ \Biggl( e^{E(P,U) + O(1)} \Biggr)^{-1} \leq \tilde{F}(h) \leq e^{E(P,U) + O(1)} . $$
Furthermore, the value of $\tilde{F}(h)$ is independent of all the products $|I_{l, \text{rand}}(\tilde{h}(l))|$ with $0 \leq l \leq \log\log P - 10\log\log\log P - 1$, and the $f(p)$ appearing in them. So conditioning on the value $F$ of $\tilde{F}(h)$, we have the upper bound
$$ \E \textbf{1}_{\tilde{\mathcal{G}}^{'}(h)} |\prod_{p \leq P} (1 - \frac{f(p)}{p^{1/2+ih}})|^{-2} \leq \sup_{|\log F| \leq E(P,U) + O(1)} \E \textbf{1}_{\tilde{\mathcal{G}}^{''}_{F}(h)} |\prod_{p \leq P} (1 - \frac{f(p)}{p^{1/2+ih}})|^{-2} , $$
where $\tilde{\mathcal{G}}^{''}_{F}(h)$ is the event that: for all $0 \leq j \leq \log\log P - 10\log\log\log P - 1$ we have
$$ F \Biggl( \frac{\log P}{e^j} e^{3\log\log\log P + U + 19/3} \Biggr)^{-1} \leq \prod_{l = j}^{\lfloor \log\log P - 10\log\log\log P \rfloor - 1} |I_{l, \text{rand}}(\tilde{h}(l))| \leq F \frac{\log P}{e^j} e^{3\log\log\log P + U + 19/3} , $$
and also
$$ \prod_{l = 0}^{\lfloor \log\log P - 10\log\log\log P \rfloor - 1} |I_{l, \text{rand}}(\tilde{h}(l))| \geq F \frac{\log P}{V B_0 e^{19/3}} . $$

Since all the primes $p$ involved in these $I_{l, \text{rand}}(\tilde{h}(l))$ satisfy the lower bound $\log p > e^{-(l+1)}\log P \geq e^{- \lfloor \log\log P - 10\log\log\log P \rfloor} \log P \geq (\log\log P)^{10}$, and in our barrier conditions in the definition of $\tilde{\mathcal{G}}^{''}_{F}(h)$ we have $3\log\log\log P + U \ll \log\log P$ and $|\log F| \ll \log\log\log P + U \ll \log\log P$, the size restrictions (relating $x_j, u_j, v_j$) in Lemma 4 of Harper~\cite{harperrmflow} are now amply satisfied. It implies that $\E \textbf{1}_{\tilde{\mathcal{G}}^{''}_{F}(h)} |\prod_{p \leq P} (1 - \frac{f(p)}{p^{1/2+ih}})|^{-2}$ is
\begin{eqnarray}
& \ll & \p\Biggl(\sum_{m=1}^{j} G_m \leq \log F + E(P,U) + O(1) \; \forall \; 1 \leq j \leq \log\log P - 10\log\log\log P, \; \text{and}      \nonumber \\
&& \sum_{m=1}^{\lfloor \log\log P - 10\log\log\log P \rfloor} G_m \geq \log F - \log V + 10\log\log\log P + O(1) \Biggr) \cdot \E |\prod_{p \leq P} (1 - \frac{f(p)}{p^{1/2+ih}})|^{-2} , \nonumber
\end{eqnarray}
where the $G_m$ are independent Gaussian random variables, each having mean 0 and variance $1/2 + o(1)$ (as $P \rightarrow \infty$). Note that in the application of Lemma 4, all the barriers are shifted by $10\log\log\log P + O(1)$ to account for the fact that we now ``start'' our Euler products at $\lfloor \log\log P - 10\log\log\log P \rfloor - 1$ rather than at $\log\log P$ (so we ``start'' with $(\log P)/e^j \asymp e^{10\log\log\log P}$ rather than $\asymp 1$). Note also that Lemma 4 of Harper~\cite{harperrmflow} actually implies a slightly stronger upper bound where the probability involves a (weak) lower bound condition on all the $\sum_{m=1}^{j} G_m$ as well, but we discard that for simplicity.

Finally, we already observed in the proof of Lemma 4 above that $\E |\prod_{p \leq P} (1 - \frac{f(p)}{p^{1/2+ih}})|^{-2} \asymp \log P$. And Probability Result 1, applied with $a = \log F + E(P,U) + O(1)$ and $b = E(P,U) + \log V - 10\log\log\log P + O(1)$ and $n \asymp \log\log P$, implies that our probability is
\begin{eqnarray}
& \ll & \min\{1, \frac{\log F + E(P,U) + O(1)}{\sqrt{\log\log P}}\} (\min\{1, \frac{E(P,U) + \log V - 10\log\log\log P}{\sqrt{\log\log P}}\})^2 \nonumber \\
& \ll & \min\{1, \frac{\log\log\log P + U}{\sqrt{\log\log P}} \} \min\{1, \frac{\log\log\log P + U + \log V}{\sqrt{\log\log P}} \}^2 , \nonumber
\end{eqnarray}
as required.
\end{proof}

\begin{proof}[Proof of Key Proposition 4]
We may assume without loss of generality that $e^{-U} \leq V \leq e^{\sqrt{\log\log P}}$, since if $V$ is larger we just discard the indicator $\textbf{1}_{|\sum_{\substack{m \leq T^{\epsilon}, \\ P \; \text{smooth}}} \frac{1}{m^{1/2+i(t+h)}}| > (\log P)/V}$ from the integral in Key Proposition 4 and proceed as in the proof of Key Proposition 1 (or, if the reader prefers, proceed as in the proof below but omitting all the new manipulations related to this indicator function).

Now we note that if the event $\tilde{\mathcal{G}}_{t}(h)$ occurs, then in particular (using the series expansion of the logarithm and the condition \eqref{tilgdef} for $j=0$) we have
$$ \Biggl|\sum_{p^{j} \leq P} \frac{1}{j p^{j(1/2 + it + i\tilde{h}(l(p)))}} \Biggr| = \Biggl|\sum_{l=0}^{\lfloor \log\log P \rfloor - B - 1} \log I_{l,t}(\tilde{h}(l)) + O(1) \Biggr| \leq 4\log\log P , $$
say. Here the $O(1)$ term depends on the fixed constant $B$. Using the series expansion of the exponential and our condition that $\frac{\epsilon \log T}{\log P} > 20\log\log P$, this implies that
\begin{eqnarray}
\sum_{k=0}^{\lfloor \frac{\epsilon \log T}{\log P} \rfloor} \frac{1}{k!} (\sum_{p^{j} \leq P} \frac{1}{j p^{j(1/2 + it + i\tilde{h}(l(p)))}} )^k & = & \exp\{ \sum_{p^{j} \leq P} \frac{1}{j p^{j(1/2 + it + i\tilde{h}(l(p)))}} \} + O\Biggl(\frac{(4\log\log P)^{\lfloor \frac{\epsilon \log T}{\log P} \rfloor + 1}}{(\lfloor \frac{\epsilon \log T}{\log P} \rfloor + 1)!} \Biggr) \nonumber \\
& = & e^{O(1)} \prod_{l = 0}^{\lfloor \log\log P \rfloor - B - 1} I_{l, t}(\tilde{h}(l)) + O(1) . \nonumber
\end{eqnarray}
This is where we benefit from the fact that \eqref{tilgdef} deals with $\log I_{l,t}(\tilde{h}(l))$ rather than just $\log |I_{l,t}(\tilde{h}(l))|$. Thus for a certain absolute constant $B_0$ (whose value could be given in terms of the constant $B$), we can upper bound the integrand in Key Proposition 4:
\begin{eqnarray}
&& \textbf{1}_{\tilde{\mathcal{G}}_{t}(h)} \textbf{1}_{|\sum_{\substack{m \leq T^{\epsilon}, \\ P \; \text{smooth}}} \frac{1}{m^{1/2+i(t+h)}}| > (\log P)/V} \Biggl|\sum_{\substack{m \leq T^{\epsilon}, \\ P \; \text{smooth}}} \frac{1}{m^{1/2+i(t+h)}} \Biggr|^2 \Biggl|\sum_{\substack{n \leq T^{1/2 - 2\epsilon}, \\ n \; \text{is} \; P \; \text{rough}}} \frac{1}{n^{\sigma+i(t+h)}} \Biggr|^2 \nonumber \\
& \ll & \textbf{1}_{\tilde{\mathcal{G}}_{t}(h)} \textbf{1}_{\prod_{l = 0}^{\lfloor \log\log P \rfloor - B - 1} |I_{l, t}(\tilde{h}(l))| > \frac{\log P}{B_0 V}} \Biggl|\sum_{\substack{m \leq T^{\epsilon}, \\ P \; \text{smooth}}} \frac{1}{m^{1/2+i(t+h)}} \Biggr|^2 \Biggl|\sum_{\substack{n \leq T^{1/2 - 2\epsilon}, \\ n \; \text{is} \; P \; \text{rough}}} \frac{1}{n^{\sigma+i(t+h)}} \Biggr|^2 + \nonumber \\
&& + \Biggl|\sum_{\substack{m \leq T^{\epsilon}, \\ P \; \text{smooth}}} \frac{1}{m^{1/2+i(t+h)}} - \sum_{k=0}^{\lfloor \frac{\epsilon \log T}{\log P} \rfloor} \frac{1}{k!} (\sum_{p^{j} \leq P} \frac{1}{j p^{j(1/2 + it + i\tilde{h}(l(p)))}} )^k \Biggr|^2 \Biggl|\sum_{\substack{n \leq T^{1/2 - 2\epsilon}, \\ n \; \text{is} \; P \; \text{rough}}} \frac{1}{n^{\sigma+i(t+h)}} \Biggr|^2 , \nonumber
\end{eqnarray}
since if $\prod_{l = 0}^{\lfloor \log\log P \rfloor - B - 1} |I_{l, t}(\tilde{h}(l))| \leq \frac{\log P}{B_0 V}$ but $|\sum_{\substack{m \leq T^{\epsilon}, \\ P \; \text{smooth}}} \frac{1}{m^{1/2+i(t+h)}}| > \frac{\log P}{V}$, then the above calculations (and our upper bound $V \leq e^{\sqrt{\log\log P}}$) imply $|\sum_{\substack{m \leq T^{\epsilon}, \\ P \; \text{smooth}}} \frac{1}{m^{1/2+i(t+h)}}|$ has comparable magnitude to the difference $|\sum_{\substack{m \leq T^{\epsilon}, \\ P \; \text{smooth}}} \frac{1}{m^{1/2+i(t+h)}} - \sum_{k=0}^{\lfloor \frac{\epsilon \log T}{\log P} \rfloor} \frac{1}{k!} (\sum_{p^{j} \leq P} \frac{1}{j p^{j(1/2 + it + i\tilde{h}(l(p)))}} )^k |$.

By Lemma \ref{eulerprodapproxlem}, the integral of the second term over $T \leq t \leq 2T$ is $\ll \frac{T \log T}{(\log\log P)^2}$, which is more than good enough for Key Proposition 4. To bound the integral of the first term, we upper bound $\textbf{1}_{\tilde{\mathcal{G}}_{t}(h)} \textbf{1}_{\prod_{l = 0}^{\lfloor \log\log P \rfloor - B - 1} |I_{l, t}(\tilde{h}(l))| > \frac{\log P}{B_0 V}}$ by a product of smooth functions. When doing this, we can afford to first replace \eqref{tilgdef} in the definition of $\tilde{\mathcal{G}}_{t}(h)$ by the weaker condition that
$$ \Biggl( \frac{\log P}{e^j} e^{3\log\log\log P + U} \Biggr)^{-1} \leq \prod_{l = j}^{\lfloor \log\log P \rfloor - B - 1} |I_{l,t}(\tilde{h}(l))| \leq \frac{\log P}{e^j} e^{3\log\log\log P + U} , $$
and so the situation is exactly similar as in the proof of Key Proposition 1 except that, when we introduce the smooth approximating function corresponding to $j=0$, we need to shift its argument and modify the choice of the parameter $R_0$ to incorporate the stronger lower bound condition $\prod_{l = 0}^{\lfloor \log\log P \rfloor - B - 1} |I_{l, t}(\tilde{h}(l))| > \frac{\log P}{B_0 V}$. Everything proceeds in the same way as for Key Proposition 1, and we finally deduce that (up to error terms of the same form as \eqref{mainerrorkp1} and \eqref{othererrorkp1}) our first integral (divided by $1/T$) is
$$ \ll (1 + \delta)^{\log\log P} \E \left(\textbf{1}_{\tilde{\mathcal{G}}^{'}(h)} + \delta \right) |\sum_{\substack{m \leq T^{\epsilon}, \\ m \; \text{is} \; P \; \text{smooth}}} \frac{f(m)}{m^{1/2+ih}}|^2 |\sum_{\substack{n \leq T^{1/2 - 2\epsilon}, \\ n \; \text{is} \; P \; \text{rough}}} \frac{f(n)}{n^{\sigma+ih}}|^2 . $$
Here $\delta > 0$ is a small parameter that we may choose, and we must satisfy the condition $T^{1/2 - \epsilon} P^{2k(\lfloor \log\log P \rfloor + 1)} < T$, where $k = \lfloor (\frac{100\log\log P}{\delta})^2 \rfloor$.

Finally, if we choose $\delta = 1/(\log\log P)^2$ (say) then $k \asymp (\log\log P)^6$, so the condition $T^{1/2 - \epsilon} P^{2k(\lfloor \log\log P \rfloor + 1)} < T$ is amply satisfied given our assumption that $P \leq T^{1/(\log\log T)^8}$, and the error terms \eqref{mainerrorkp1}, \eqref{othererrorkp1} are very small. Furthermore, as in the proof of Key Proposition 1 we have $(1 + \delta)^{\log\log P} \ll 1$ and $\delta \E |\sum_{\substack{m \leq T^{\epsilon}, \\ P \; \text{smooth}}} \frac{f(m)}{m^{1/2+ih}}|^2 |\sum_{\substack{n \leq T^{1/2 - 2\epsilon}, \\ P \; \text{rough}}} \frac{f(n)}{n^{\sigma+ih}}|^2 \asymp \delta \log T = \frac{\log T}{(\log\log P)^2}$, which is more than good enough for Key Proposition 4. And the contribution from the term involving $\textbf{1}_{\tilde{\mathcal{G}}^{'}(h)}$ may be estimated satisfactorily using Lemma \ref{mainprob2lem}, together with our condition that $\epsilon > \frac{20\log P \log\log P}{\log T}$.
\end{proof}

\vspace{12pt}
\noindent {\em Acknowledgements.}
The author would like to thank Louis-Pierre Arguin, for informing him of his forthcoming joint work on the Fyodorov--Hiary--Keating conjecture. He would also like to thank Jon Keating, Joseph Najnudel and K. Soundararajan for their comments and encouragement.

\end{document}